\documentclass[12pt]{amsart}
\title[Markov elements]{Markov elements in   affine Temperley-Lieb algebras}
\author{Sadek AL HARBAT}
\address{LAMFA, Universit\'e de Picardie - Jules Verne} 
\email{sadikharbat@math.univ-paris-diderot.fr}

\usepackage{etex} 
\usepackage[latin1]{inputenc}	
\usepackage[T1]{fontenc}
\usepackage[english]{babel}
\usepackage{amsmath,amssymb}
\usepackage{wasysym}
\usepackage{stmaryrd}
\usepackage[table]{xcolor}
\usepackage{color}
\usepackage{dsfont}\let\mathbb\mathds
\usepackage[all]{xy}
\usepackage{graphicx}
\usepackage{mathenv}
\usepackage{lastpage}
\usepackage{fancyhdr}	
\usepackage{subfigure}
\usepackage{geometry}
\usepackage[version=3]{mhchem}
\usepackage[force,almostfull]{textcomp}
\usepackage{array}
\usepackage{lipsum}		
\usepackage{eso-pic}
\usepackage{multirow}
\usepackage{fancybox}
\usepackage{cite}
\usepackage{amsthm}
\usepackage{listings}
\usepackage{lscape}
\usepackage{multicol}		
\usepackage{setspace}
\usepackage{times}
\usepackage{eurosym}
\usepackage{latexsym}
\usepackage{ulem}
\usepackage{lmodern}
\usepackage{colortbl}
\usepackage{rotating}
\usepackage{type1cm}
\usepackage{lettrine}
\usepackage{ragged2e}
\usepackage{mathrsfs}
\usepackage{scrtime}
\usepackage{enumitem}

\usepackage{longtable}
\usepackage{url}
\usepackage{nomentbl}
\usepackage{nomencl}
\usepackage{hyperref}
\hypersetup{
pdfpagemode=UseOutlines,      
pdfstartview=Fit,             
pdffitwindow=true,            
pdfpagelayout=TwoColumnsRight,
pdftoolbar=true,              
pdfmenubar=true,              
bookmarksopen=false,          
bookmarksnumbered=true,       
colorlinks=true,              
pdfauthor={Ton nom},          
pdftitle={Titre PDF},         
pdfcreator=PDFLaTeX,          %
pdfproducer=PDFLaTeX,         %
linkcolor=blue,               
urlcolor=blue,                
anchorcolor=black,            
citecolor=blue,               
frenchlinks=true,             
pdfborder={0 0 0}             
}
\usepackage{pgf, tikz, tikz-cd}
\usetikzlibrary{matrix,positioning}
\usetikzlibrary{arrows,decorations.markings}

\definecolor{aliceblue}{rgb}{0.94,0.97,1.00}
\definecolor{antiquewhite}{rgb}{0.98,0.92,0.84}
\definecolor{antiquewhite1}{rgb}{1.00,0.94,0.86}
\definecolor{antiquewhite2}{rgb}{0.93,0.87,0.80}
\definecolor{antiquewhite3}{rgb}{0.80,0.75,0.69}
\definecolor{antiquewhite4}{rgb}{0.55,0.51,0.47}
\definecolor{aquamarine}{rgb}{0.50,1.00,0.83}
\definecolor{aquamarine1}{rgb}{0.50,1.00,0.83}
\definecolor{aquamarine2}{rgb}{0.46,0.93,0.78}
\definecolor{aquamarine3}{rgb}{0.40,0.80,0.67}
\definecolor{aquamarine4}{rgb}{0.27,0.55,0.45}
\definecolor{azure}{rgb}{0.94,1.00,1.00}
\definecolor{azure1}{rgb}{0.94,1.00,1.00}
\definecolor{azure2}{rgb}{0.88,0.93,0.93}
\definecolor{azure3}{rgb}{0.76,0.80,0.80}
\definecolor{azure4}{rgb}{0.51,0.55,0.55}
\definecolor{beige}{rgb}{0.96,0.96,0.86}
\definecolor{bisque}{rgb}{1.00,0.89,0.77}
\definecolor{bisque1}{rgb}{1.00,0.89,0.77}
\definecolor{bisque2}{rgb}{0.93,0.84,0.72}
\definecolor{bisque3}{rgb}{0.80,0.72,0.62}
\definecolor{bisque4}{rgb}{0.55,0.49,0.42}
\definecolor{black}{rgb}{0.00,0.00,0.00}
\definecolor{blanchedalmond}{rgb}{1.00,0.92,0.80}
\definecolor{blue}{rgb}{0.00,0.00,1.00}
\definecolor{blue1}{rgb}{0.00,0.00,1.00}
\definecolor{blue2}{rgb}{0.00,0.00,0.93}
\definecolor{blue3}{rgb}{0.00,0.00,0.80}
\definecolor{blue4}{rgb}{0.00,0.00,0.55}
\definecolor{blueviolet}{rgb}{0.54,0.17,0.89}
\definecolor{brown}{rgb}{0.65,0.16,0.16}
\definecolor{brown1}{rgb}{1.00,0.25,0.25}
\definecolor{brown2}{rgb}{0.93,0.23,0.23}
\definecolor{brown3}{rgb}{0.80,0.20,0.20}
\definecolor{brown4}{rgb}{0.55,0.14,0.14}
\definecolor{burlywood}{rgb}{0.87,0.72,0.53}
\definecolor{burlywood1}{rgb}{1.00,0.83,0.61}
\definecolor{burlywood2}{rgb}{0.93,0.77,0.57}
\definecolor{burlywood3}{rgb}{0.80,0.67,0.49}
\definecolor{burlywood4}{rgb}{0.55,0.45,0.33}
\definecolor{cadetblue}{rgb}{0.37,0.62,0.63}
\definecolor{cadetblue1}{rgb}{0.60,0.96,1.00}
\definecolor{cadetblue2}{rgb}{0.56,0.90,0.93}
\definecolor{cadetblue3}{rgb}{0.48,0.77,0.80}
\definecolor{cadetblue4}{rgb}{0.33,0.53,0.55}
\definecolor{chartreuse}{rgb}{0.50,1.00,0.00}
\definecolor{chartreuse1}{rgb}{0.50,1.00,0.00}
\definecolor{chartreuse2}{rgb}{0.46,0.93,0.00}
\definecolor{chartreuse3}{rgb}{0.40,0.80,0.00}
\definecolor{chartreuse4}{rgb}{0.27,0.55,0.00}
\definecolor{chocolate}{rgb}{0.82,0.41,0.12}
\definecolor{chocolate1}{rgb}{1.00,0.50,0.14}
\definecolor{chocolate2}{rgb}{0.93,0.46,0.13}
\definecolor{chocolate3}{rgb}{0.80,0.40,0.11}
\definecolor{chocolate4}{rgb}{0.55,0.27,0.07}
\definecolor{coral}{rgb}{1.00,0.50,0.31}
\definecolor{coral1}{rgb}{1.00,0.45,0.34}
\definecolor{coral2}{rgb}{0.93,0.42,0.31}
\definecolor{coral3}{rgb}{0.80,0.36,0.27}
\definecolor{coral4}{rgb}{0.55,0.24,0.18}
\definecolor{cornsilk}{rgb}{1.00,0.97,0.86}
\definecolor{cornflowerblue}{rgb}{0.39,0.58,0.93}
\definecolor{cornsilk1}{rgb}{1.00,0.97,0.86}
\definecolor{cornsilk2}{rgb}{0.93,0.91,0.80}
\definecolor{cornsilk3}{rgb}{0.80,0.78,0.69}
\definecolor{cornsilk4}{rgb}{0.55,0.53,0.47}
\definecolor{cyan}{rgb}{0.00,1.00,1.00}
\definecolor{cyan1}{rgb}{0.00,1.00,1.00}
\definecolor{cyan2}{rgb}{0.00,0.93,0.93}
\definecolor{cyan3}{rgb}{0.00,0.80,0.80}
\definecolor{cyan4}{rgb}{0.00,0.55,0.55}
\definecolor{darkblue}{rgb}{0.00,0.00,0.55}
\definecolor{darkcyan}{rgb}{0.00,0.55,0.55}
\definecolor{darkgoldenrod}{rgb}{0.72,0.53,0.04}
\definecolor{DarkGoldenrod1}{rgb}{1.00,0.73,0.06}
\definecolor{DarkGoldenrod2}{rgb}{0.93,0.68,0.05}
\definecolor{DarkGoldenrod3}{rgb}{0.80,0.58,0.05}
\definecolor{DarkGoldenrod4}{rgb}{0.55,0.40,0.03}
\definecolor{darkgray}{rgb}{0.66,0.66,0.66}
\definecolor{darkgreen}{rgb}{0.00,0.39,0.00}
\definecolor{darkgrey}{rgb}{0.66,0.66,0.66}
\definecolor{darkkhaki}{rgb}{0.74,0.72,0.42}
\definecolor{darkmagenta}{rgb}{0.55,0.00,0.55}
\definecolor{darkOlivegreen}{rgb}{0.33,0.42,0.18}
\definecolor{darkOlivegreen1}{rgb}{0.79,1.00,0.44}
\definecolor{darkOlivegreen2}{rgb}{0.74,0.93,0.41}
\definecolor{darkOlivegreen3}{rgb}{0.64,0.80,0.35}
\definecolor{darkOlivegreen4}{rgb}{0.43,0.55,0.24}
\definecolor{darkolive}{rgb}{0.33,0.42,0.18}
\definecolor{DarkOrange1}{rgb}{1.00,0.50,0.00}
\definecolor{DarkOrange2}{rgb}{0.93,0.46,0.00}
\definecolor{DarkOrange3}{rgb}{0.80,0.40,0.00}
\definecolor{DarkOrange4}{rgb}{0.55,0.27,0.00}
\definecolor{darkorange}{rgb}{1.00,0.55,0.00}
\definecolor{darkorchid}{rgb}{0.60,0.20,0.80}
\definecolor{darkOrchid1}{rgb}{0.75,0.24,1.00}
\definecolor{darkOrchid2}{rgb}{0.70,0.23,0.93}
\definecolor{darkOrchid3}{rgb}{0.60,0.20,0.80}
\definecolor{darkOrchid4}{rgb}{0.41,0.13,0.55}
\definecolor{darkred}{rgb}{0.55,0.00,0.00}
\definecolor{darksalmon}{rgb}{0.91,0.59,0.48}
\definecolor{darkseagreen}{rgb}{0.56,0.74,0.56}
\definecolor{darkseagreen1}{rgb}{0.76,1.00,0.76}
\definecolor{darkseagreen2}{rgb}{0.71,0.93,0.71}
\definecolor{darkseagreen3}{rgb}{0.61,0.80,0.61}
\definecolor{darkseagreen4}{rgb}{0.41,0.55,0.41}
\definecolor{darksea}{rgb}{0.56,0.74,0.56}
\definecolor{darkslategray}{rgb}{0.18,0.31,0.31}
\definecolor{darkslateblue}{rgb}{0.28,0.24,0.55}
\definecolor{darkslategray1}{rgb}{0.59,1.00,1.00}
\definecolor{darkslategray2}{rgb}{0.55,0.93,0.93}
\definecolor{darkslategray3}{rgb}{0.47,0.80,0.80}
\definecolor{darkslategray4}{rgb}{0.32,0.55,0.55}
\definecolor{darkslate}{rgb}{0.18,0.31,0.31}
\definecolor{darkslate1}{rgb}{0.28,0.24,0.55}
\definecolor{darkturquoise}{rgb}{0.00,0.81,0.82}
\definecolor{darkviolet}{rgb}{0.58,0.00,0.83}
\definecolor{deeppink}{rgb}{1.00,0.08,0.58}
\definecolor{deepPink1}{rgb}{1.00,0.08,0.58}
\definecolor{deepPink2}{rgb}{0.93,0.07,0.54}
\definecolor{deepPink3}{rgb}{0.80,0.06,0.46}
\definecolor{deepPink4}{rgb}{0.55,0.04,0.31}
\definecolor{deepskyblue}{rgb}{0.00,0.75,1.00}
\definecolor{deepskyblue1}{rgb}{0.00,0.75,1.00}
\definecolor{deepskyblue2}{rgb}{0.00,0.70,0.93}
\definecolor{deepskyblue3}{rgb}{0.00,0.60,0.80}
\definecolor{deepskyblue4}{rgb}{0.00,0.41,0.55}
\definecolor{deepsky}{rgb}{0.00,0.75,1.00}
\definecolor{dimgray}{rgb}{0.41,0.41,0.41}
\definecolor{dodgerblue}{rgb}{0.12,0.56,1.00}
\definecolor{dodgerblue1}{rgb}{0.12,0.56,1.00}
\definecolor{dodgerblue2}{rgb}{0.11,0.53,0.93}
\definecolor{dodgerblue3}{rgb}{0.09,0.45,0.80}
\definecolor{dodgerblue4}{rgb}{0.06,0.31,0.55}
\definecolor{firebrick}{rgb}{0.70,0.13,0.13}
\definecolor{firebrick1}{rgb}{1.00,0.19,0.19}
\definecolor{firebrick2}{rgb}{0.93,0.17,0.17}
\definecolor{firebrick3}{rgb}{0.80,0.15,0.15}
\definecolor{firebrick4}{rgb}{0.55,0.10,0.10}
\definecolor{floralwhite}{rgb}{1.00,0.98,0.94}
\definecolor{forestgreen}{rgb}{0.13,0.55,0.13}
\definecolor{gainsboro}{rgb}{0.86,0.86,0.86}
\definecolor{ghostwhite}{rgb}{0.97,0.97,1.00}
\definecolor{gold1}{rgb}{1.00,0.84,0.00}
\definecolor{gold2}{rgb}{0.93,0.79,0.00}
\definecolor{gold3}{rgb}{0.80,0.68,0.00}
\definecolor{gold4}{rgb}{0.55,0.46,0.00}
\definecolor{goldenrod}{rgb}{0.85,0.65,0.13}
\definecolor{goldenrod1}{rgb}{1.00,0.76,0.15}
\definecolor{goldenrod2}{rgb}{0.93,0.71,0.13}
\definecolor{goldenrod3}{rgb}{0.80,0.61,0.11}
\definecolor{goldenrod4}{rgb}{0.55,0.41,0.08}
\definecolor{gold}{rgb}{1.00,0.84,0.00}
\definecolor{gray}{rgb}{0.75,0.75,0.75}
\definecolor{gray0}{rgb}{0.00,0.00,0.00}
\definecolor{gray1}{rgb}{0.01,0.01,0.01}
\definecolor{gray2}{rgb}{0.02,0.02,0.02}
\definecolor{gray3}{rgb}{0.03,0.03,0.03}
\definecolor{gray4}{rgb}{0.04,0.04,0.04}
\definecolor{gray5}{rgb}{0.05,0.05,0.05}
\definecolor{gray6}{rgb}{0.06,0.06,0.06}
\definecolor{gray7}{rgb}{0.07,0.07,0.07}
\definecolor{gray8}{rgb}{0.08,0.08,0.08}
\definecolor{gray9}{rgb}{0.09,0.09,0.09}
\definecolor{gray10}{rgb}{0.10,0.10,0.10}
\definecolor{gray11}{rgb}{0.11,0.11,0.11}
\definecolor{gray12}{rgb}{0.12,0.12,0.12}
\definecolor{gray13}{rgb}{0.13,0.13,0.13}
\definecolor{gray14}{rgb}{0.14,0.14,0.14}
\definecolor{gray15}{rgb}{0.15,0.15,0.15}
\definecolor{gray16}{rgb}{0.16,0.16,0.16}
\definecolor{gray17}{rgb}{0.17,0.17,0.17}
\definecolor{gray18}{rgb}{0.18,0.18,0.18}
\definecolor{gray19}{rgb}{0.19,0.19,0.19}
\definecolor{gray20}{rgb}{0.20,0.20,0.20}
\definecolor{gray21}{rgb}{0.21,0.21,0.21}
\definecolor{gray22}{rgb}{0.22,0.22,0.22}
\definecolor{gray23}{rgb}{0.23,0.23,0.23}
\definecolor{gray24}{rgb}{0.24,0.24,0.24}
\definecolor{gray25}{rgb}{0.25,0.25,0.25}
\definecolor{gray26}{rgb}{0.26,0.26,0.26}
\definecolor{gray27}{rgb}{0.27,0.27,0.27}
\definecolor{gray28}{rgb}{0.28,0.28,0.28}
\definecolor{gray29}{rgb}{0.29,0.29,0.29}
\definecolor{gray30}{rgb}{0.30,0.30,0.30}
\definecolor{gray31}{rgb}{0.31,0.31,0.31}
\definecolor{gray32}{rgb}{0.32,0.32,0.32}
\definecolor{gray33}{rgb}{0.33,0.33,0.33}
\definecolor{gray34}{rgb}{0.34,0.34,0.34}
\definecolor{gray35}{rgb}{0.35,0.35,0.35}
\definecolor{gray36}{rgb}{0.36,0.36,0.36}
\definecolor{gray37}{rgb}{0.37,0.37,0.37}
\definecolor{gray38}{rgb}{0.38,0.38,0.38}
\definecolor{gray39}{rgb}{0.39,0.39,0.39}
\definecolor{gray40}{rgb}{0.40,0.40,0.40}
\definecolor{gray41}{rgb}{0.41,0.41,0.41}
\definecolor{gray42}{rgb}{0.42,0.42,0.42}
\definecolor{gray43}{rgb}{0.43,0.43,0.43}
\definecolor{gray44}{rgb}{0.44,0.44,0.44}
\definecolor{gray45}{rgb}{0.45,0.45,0.45}
\definecolor{gray46}{rgb}{0.46,0.46,0.46}
\definecolor{gray47}{rgb}{0.47,0.47,0.47}
\definecolor{gray48}{rgb}{0.48,0.48,0.48}
\definecolor{gray49}{rgb}{0.49,0.49,0.49}
\definecolor{gray50}{rgb}{0.50,0.50,0.50}
\definecolor{gray51}{rgb}{0.51,0.51,0.51}
\definecolor{gray52}{rgb}{0.52,0.52,0.52}
\definecolor{gray53}{rgb}{0.53,0.53,0.53}
\definecolor{gray54}{rgb}{0.54,0.54,0.54}
\definecolor{gray55}{rgb}{0.55,0.55,0.55}
\definecolor{gray56}{rgb}{0.56,0.56,0.56}
\definecolor{gray57}{rgb}{0.57,0.57,0.57}
\definecolor{gray58}{rgb}{0.58,0.58,0.58}
\definecolor{gray59}{rgb}{0.59,0.59,0.59}
\definecolor{gray60}{rgb}{0.60,0.60,0.60}
\definecolor{gray61}{rgb}{0.61,0.61,0.61}
\definecolor{gray62}{rgb}{0.62,0.62,0.62}
\definecolor{gray63}{rgb}{0.63,0.63,0.63}
\definecolor{gray64}{rgb}{0.64,0.64,0.64}
\definecolor{gray65}{rgb}{0.65,0.65,0.65}
\definecolor{gray66}{rgb}{0.66,0.66,0.66}
\definecolor{gray67}{rgb}{0.67,0.67,0.67}
\definecolor{gray68}{rgb}{0.68,0.68,0.68}
\definecolor{gray69}{rgb}{0.69,0.69,0.69}
\definecolor{gray70}{rgb}{0.70,0.70,0.70}
\definecolor{gray71}{rgb}{0.71,0.71,0.71}
\definecolor{gray72}{rgb}{0.72,0.72,0.72}
\definecolor{gray73}{rgb}{0.73,0.73,0.73}
\definecolor{gray74}{rgb}{0.74,0.74,0.74}
\definecolor{gray75}{rgb}{0.75,0.75,0.75}
\definecolor{gray76}{rgb}{0.76,0.76,0.76}
\definecolor{gray77}{rgb}{0.77,0.77,0.77}
\definecolor{gray78}{rgb}{0.78,0.78,0.78}
\definecolor{gray79}{rgb}{0.79,0.79,0.79}
\definecolor{gray80}{rgb}{0.80,0.80,0.80}
\definecolor{gray81}{rgb}{0.81,0.81,0.81}
\definecolor{gray82}{rgb}{0.82,0.82,0.82}
\definecolor{gray83}{rgb}{0.83,0.83,0.83}
\definecolor{gray84}{rgb}{0.84,0.84,0.84}
\definecolor{gray85}{rgb}{0.85,0.85,0.85}
\definecolor{gray86}{rgb}{0.86,0.86,0.86}
\definecolor{gray87}{rgb}{0.87,0.87,0.87}
\definecolor{gray88}{rgb}{0.88,0.88,0.88}
\definecolor{gray89}{rgb}{0.89,0.89,0.89}
\definecolor{gray90}{rgb}{0.90,0.90,0.90}
\definecolor{gray91}{rgb}{0.91,0.91,0.91}
\definecolor{gray92}{rgb}{0.92,0.92,0.92}
\definecolor{gray93}{rgb}{0.93,0.93,0.93}
\definecolor{gray94}{rgb}{0.94,0.94,0.94}
\definecolor{gray95}{rgb}{0.95,0.95,0.95}
\definecolor{gray96}{rgb}{0.96,0.96,0.96}
\definecolor{gray97}{rgb}{0.97,0.97,0.97}
\definecolor{gray98}{rgb}{0.98,0.98,0.98}
\definecolor{gray99}{rgb}{0.99,0.99,0.99}
\definecolor{gray100}{rgb}{1.00,1.00,1.00}
\definecolor{green}{rgb}{0.00,1.00,0.00}
\definecolor{green1}{rgb}{0.00,1.00,0.00}
\definecolor{green2}{rgb}{0.00,0.93,0.00}
\definecolor{green3}{rgb}{0.00,0.80,0.00}
\definecolor{green4}{rgb}{0.00,0.55,0.00}
\definecolor{greenyellow}{rgb}{0.68,1.00,0.18}
\definecolor{grey}{rgb}{0.75,0.75,0.75}
\definecolor{grey0}{rgb}{0.00,0.00,0.00}
\definecolor{grey1}{rgb}{0.01,0.01,0.01}
\definecolor{grey2}{rgb}{0.02,0.02,0.02}
\definecolor{grey3}{rgb}{0.03,0.03,0.03}
\definecolor{grey4}{rgb}{0.04,0.04,0.04}
\definecolor{grey5}{rgb}{0.05,0.05,0.05}
\definecolor{grey6}{rgb}{0.06,0.06,0.06}
\definecolor{grey7}{rgb}{0.07,0.07,0.07}
\definecolor{grey8}{rgb}{0.08,0.08,0.08}
\definecolor{grey9}{rgb}{0.09,0.09,0.09}
\definecolor{grey10}{rgb}{0.10,0.10,0.10}
\definecolor{grey11}{rgb}{0.11,0.11,0.11}
\definecolor{grey12}{rgb}{0.12,0.12,0.12}
\definecolor{grey13}{rgb}{0.13,0.13,0.13}
\definecolor{grey14}{rgb}{0.14,0.14,0.14}
\definecolor{grey15}{rgb}{0.15,0.15,0.15}
\definecolor{grey16}{rgb}{0.16,0.16,0.16}
\definecolor{grey17}{rgb}{0.17,0.17,0.17}
\definecolor{grey18}{rgb}{0.18,0.18,0.18}
\definecolor{grey19}{rgb}{0.19,0.19,0.19}
\definecolor{grey20}{rgb}{0.20,0.20,0.20}
\definecolor{grey21}{rgb}{0.21,0.21,0.21}
\definecolor{grey22}{rgb}{0.22,0.22,0.22}
\definecolor{grey23}{rgb}{0.23,0.23,0.23}
\definecolor{grey24}{rgb}{0.24,0.24,0.24}
\definecolor{grey25}{rgb}{0.25,0.25,0.25}
\definecolor{grey26}{rgb}{0.26,0.26,0.26}
\definecolor{grey27}{rgb}{0.27,0.27,0.27}
\definecolor{grey28}{rgb}{0.28,0.28,0.28}
\definecolor{grey29}{rgb}{0.29,0.29,0.29}
\definecolor{grey30}{rgb}{0.30,0.30,0.30}
\definecolor{grey31}{rgb}{0.31,0.31,0.31}
\definecolor{grey32}{rgb}{0.32,0.32,0.32}
\definecolor{grey33}{rgb}{0.33,0.33,0.33}
\definecolor{grey34}{rgb}{0.34,0.34,0.34}
\definecolor{grey35}{rgb}{0.35,0.35,0.35}
\definecolor{grey36}{rgb}{0.36,0.36,0.36}
\definecolor{grey37}{rgb}{0.37,0.37,0.37}
\definecolor{grey38}{rgb}{0.38,0.38,0.38}
\definecolor{grey39}{rgb}{0.39,0.39,0.39}
\definecolor{grey40}{rgb}{0.40,0.40,0.40}
\definecolor{grey41}{rgb}{0.41,0.41,0.41}
\definecolor{grey42}{rgb}{0.42,0.42,0.42}
\definecolor{grey43}{rgb}{0.43,0.43,0.43}
\definecolor{grey44}{rgb}{0.44,0.44,0.44}
\definecolor{grey45}{rgb}{0.45,0.45,0.45}
\definecolor{grey46}{rgb}{0.46,0.46,0.46}
\definecolor{grey47}{rgb}{0.47,0.47,0.47}
\definecolor{grey48}{rgb}{0.48,0.48,0.48}
\definecolor{grey49}{rgb}{0.49,0.49,0.49}
\definecolor{grey50}{rgb}{0.50,0.50,0.50}
\definecolor{grey51}{rgb}{0.51,0.51,0.51}
\definecolor{grey52}{rgb}{0.52,0.52,0.52}
\definecolor{grey53}{rgb}{0.53,0.53,0.53}
\definecolor{grey54}{rgb}{0.54,0.54,0.54}
\definecolor{grey55}{rgb}{0.55,0.55,0.55}
\definecolor{grey56}{rgb}{0.56,0.56,0.56}
\definecolor{grey57}{rgb}{0.57,0.57,0.57}
\definecolor{grey58}{rgb}{0.58,0.58,0.58}
\definecolor{grey59}{rgb}{0.59,0.59,0.59}
\definecolor{grey60}{rgb}{0.60,0.60,0.60}
\definecolor{grey61}{rgb}{0.61,0.61,0.61}
\definecolor{grey62}{rgb}{0.62,0.62,0.62}
\definecolor{grey63}{rgb}{0.63,0.63,0.63}
\definecolor{grey64}{rgb}{0.64,0.64,0.64}
\definecolor{grey65}{rgb}{0.65,0.65,0.65}
\definecolor{grey66}{rgb}{0.66,0.66,0.66}
\definecolor{grey67}{rgb}{0.67,0.67,0.67}
\definecolor{grey68}{rgb}{0.68,0.68,0.68}
\definecolor{grey69}{rgb}{0.69,0.69,0.69}
\definecolor{grey70}{rgb}{0.70,0.70,0.70}
\definecolor{grey71}{rgb}{0.71,0.71,0.71}
\definecolor{grey72}{rgb}{0.72,0.72,0.72}
\definecolor{grey73}{rgb}{0.73,0.73,0.73}
\definecolor{grey74}{rgb}{0.74,0.74,0.74}
\definecolor{grey75}{rgb}{0.75,0.75,0.75}
\definecolor{grey76}{rgb}{0.76,0.76,0.76}
\definecolor{grey77}{rgb}{0.77,0.77,0.77}
\definecolor{grey78}{rgb}{0.78,0.78,0.78}
\definecolor{grey79}{rgb}{0.79,0.79,0.79}
\definecolor{grey80}{rgb}{0.80,0.80,0.80}
\definecolor{grey81}{rgb}{0.81,0.81,0.81}
\definecolor{grey82}{rgb}{0.82,0.82,0.82}
\definecolor{grey83}{rgb}{0.83,0.83,0.83}
\definecolor{grey84}{rgb}{0.84,0.84,0.84}
\definecolor{grey85}{rgb}{0.85,0.85,0.85}
\definecolor{grey86}{rgb}{0.86,0.86,0.86}
\definecolor{grey87}{rgb}{0.87,0.87,0.87}
\definecolor{grey88}{rgb}{0.88,0.88,0.88}
\definecolor{grey89}{rgb}{0.89,0.89,0.89}
\definecolor{grey90}{rgb}{0.90,0.90,0.90}
\definecolor{grey91}{rgb}{0.91,0.91,0.91}
\definecolor{grey92}{rgb}{0.92,0.92,0.92}
\definecolor{grey93}{rgb}{0.93,0.93,0.93}
\definecolor{grey94}{rgb}{0.94,0.94,0.94}
\definecolor{grey95}{rgb}{0.95,0.95,0.95}
\definecolor{grey96}{rgb}{0.96,0.96,0.96}
\definecolor{grey97}{rgb}{0.97,0.97,0.97}
\definecolor{grey98}{rgb}{0.98,0.98,0.98}
\definecolor{grey99}{rgb}{0.99,0.99,0.99}
\definecolor{grey100}{rgb}{1.00,1.00,1.00}
\definecolor{honeydew}{rgb}{0.94,1.00,0.94}
\definecolor{honeydew1}{rgb}{0.94,1.00,0.94}
\definecolor{honeydew2}{rgb}{0.88,0.93,0.88}
\definecolor{honeydew3}{rgb}{0.76,0.80,0.76}
\definecolor{honeydew4}{rgb}{0.51,0.55,0.51}
\definecolor{hotpink}{rgb}{1.00,0.41,0.71}
\definecolor{hotPink1}{rgb}{1.00,0.43,0.71}
\definecolor{hotPink2}{rgb}{0.93,0.42,0.65}
\definecolor{hotPink3}{rgb}{0.80,0.38,0.56}
\definecolor{hotPink4}{rgb}{0.55,0.23,0.38}
\definecolor{indianred}{rgb}{0.80,0.36,0.36}
\definecolor{indianred1}{rgb}{1.00,0.42,0.42}
\definecolor{indianred2}{rgb}{0.93,0.39,0.39}
\definecolor{indianred3}{rgb}{0.80,0.33,0.33}
\definecolor{indianred4}{rgb}{0.55,0.23,0.23}
\definecolor{ivory}{rgb}{1.00,1.00,0.94}
\definecolor{ivory1}{rgb}{1.00,1.00,0.94}
\definecolor{ivory2}{rgb}{0.93,0.93,0.88}
\definecolor{ivory3}{rgb}{0.80,0.80,0.76}
\definecolor{ivory4}{rgb}{0.55,0.55,0.51}
\definecolor{khaki}{rgb}{0.94,0.90,0.55}
\definecolor{khaki1}{rgb}{1.00,0.96,0.56}
\definecolor{khaki2}{rgb}{0.93,0.90,0.52}
\definecolor{khaki3}{rgb}{0.80,0.78,0.45}
\definecolor{khaki4}{rgb}{0.55,0.53,0.31}
\definecolor{lavenderblush}{rgb}{1.00,0.94,0.96}
\definecolor{lavenderblush1}{rgb}{1.00,0.94,0.96}
\definecolor{lavenderblush2}{rgb}{0.93,0.88,0.90}
\definecolor{lavenderblush3}{rgb}{0.80,0.76,0.77}
\definecolor{lavenderblush4}{rgb}{0.55,0.51,0.53}
\definecolor{lavender}{rgb}{0.90,0.90,0.98}
\definecolor{lawngreen}{rgb}{0.49,0.99,0.00}
\definecolor{lemonchiffon}{rgb}{1.00,0.98,0.80}
\definecolor{lemonchiffon1}{rgb}{1.00,0.98,0.80}
\definecolor{lemonchiffon2}{rgb}{0.93,0.91,0.75}
\definecolor{lemonchiffon3}{rgb}{0.80,0.79,0.65}
\definecolor{lemonchiffon4}{rgb}{0.55,0.54,0.44}
\definecolor{lightblue}{rgb}{0.68,0.85,0.90}
\definecolor{lightblue1}{rgb}{0.75,0.94,1.00}
\definecolor{lightblue2}{rgb}{0.70,0.87,0.93}
\definecolor{lightblue3}{rgb}{0.60,0.75,0.80}
\definecolor{lightblue4}{rgb}{0.41,0.51,0.55}
\definecolor{lightcoral}{rgb}{0.94,0.50,0.50}
\definecolor{lightcyan}{rgb}{0.88,1.00,1.00}
\definecolor{lightcyan1}{rgb}{0.88,1.00,1.00}
\definecolor{lightcyan2}{rgb}{0.82,0.93,0.93}
\definecolor{lightcyan3}{rgb}{0.71,0.80,0.80}
\definecolor{lightcyan4}{rgb}{0.48,0.55,0.55}
\definecolor{lightgoldenrod}{rgb}{0.93,0.87,0.51}
\definecolor{lightgoldenrod0}{rgb}{0.98,0.98,0.82}
\definecolor{lightgoldenrod1}{rgb}{1.00,0.93,0.55}
\definecolor{lightgoldenrod2}{rgb}{0.93,0.86,0.51}
\definecolor{lightgoldenrod3}{rgb}{0.80,0.75,0.44}
\definecolor{lightgoldenrod4}{rgb}{0.55,0.51,0.30}
\definecolor{lightgoldenrodYellow}{rgb}{0.98,0.98,0.82}
\definecolor{lightgray}{rgb}{0.83,0.83,0.83}
\definecolor{lightgreen}{rgb}{0.56,0.93,0.56}
\definecolor{lightgrey}{rgb}{0.83,0.83,0.83}
\definecolor{lightpink}{rgb}{1.00,0.71,0.76}
\definecolor{lightpink1}{rgb}{1.00,0.68,0.73}
\definecolor{lightpink2}{rgb}{0.93,0.64,0.68}
\definecolor{lightpink3}{rgb}{0.80,0.55,0.58}
\definecolor{lightpink4}{rgb}{0.55,0.37,0.40}
\definecolor{lightsalmon}{rgb}{1.00,0.63,0.48}
\definecolor{lightsalmon1}{rgb}{1.00,0.63,0.48}
\definecolor{lightsalmon2}{rgb}{0.93,0.58,0.45}
\definecolor{lightsalmon3}{rgb}{0.80,0.51,0.38}
\definecolor{lightsalmon4}{rgb}{0.55,0.34,0.26}
\definecolor{lightseagreen}{rgb}{0.13,0.70,0.67}
\definecolor{lightsea}{rgb}{0.13,0.70,0.67}
\definecolor{lightsky}{rgb}{0.53,0.81,0.98}
\definecolor{lightSkyblue}{rgb}{0.53,0.81,0.98}
\definecolor{lightSkyblue1}{rgb}{0.69,0.89,1.00}
\definecolor{lightSkyblue2}{rgb}{0.64,0.83,0.93}
\definecolor{lightSkyblue3}{rgb}{0.55,0.71,0.80}
\definecolor{lightSkyblue4}{rgb}{0.38,0.48,0.55}
\definecolor{lightslateblue}{rgb}{0.52,0.44,1.00}
\definecolor{lightslategray}{rgb}{0.47,0.53,0.60}
\definecolor{lightslate}{rgb}{0.47,0.53,0.60}
\definecolor{lightslate1}{rgb}{0.52,0.44,1.00}
\definecolor{lightsteelblue}{rgb}{0.69,0.77,0.87}
\definecolor{lightsteelblue1}{rgb}{0.79,0.88,1.00}
\definecolor{lightsteelblue2}{rgb}{0.74,0.82,0.93}
\definecolor{lightsteelblue3}{rgb}{0.64,0.71,0.80}
\definecolor{lightsteelblue4}{rgb}{0.43,0.48,0.55}
\definecolor{lightyellow}{rgb}{1.00,1.00,0.88}
\definecolor{lightsteel}{rgb}{0.69,0.77,0.87}
\definecolor{lightyellow}{rgb}{1.00,1.00,0.88}
\definecolor{lightyellow1}{rgb}{1.00,1.00,0.88}
\definecolor{lightyellow2}{rgb}{0.93,0.93,0.82}
\definecolor{lightyellow3}{rgb}{0.80,0.80,0.71}
\definecolor{Lightyellow4}{rgb}{0.55,0.55,0.48}
\definecolor{limegreen}{rgb}{0.20,0.80,0.20}
\definecolor{linen}{rgb}{0.98,0.94,0.90}
\definecolor{magenta}{rgb}{1.00,0.00,1.00}
\definecolor{magenta1}{rgb}{1.00,0.00,1.00}
\definecolor{magenta2}{rgb}{0.93,0.00,0.93}
\definecolor{magenta3}{rgb}{0.80,0.00,0.80}
\definecolor{magenta4}{rgb}{0.55,0.00,0.55}
\definecolor{maroon}{rgb}{0.69,0.19,0.38}
\definecolor{maroon1}{rgb}{1.00,0.20,0.70}
\definecolor{maroon2}{rgb}{0.93,0.19,0.65}
\definecolor{maroon3}{rgb}{0.80,0.16,0.56}
\definecolor{maroon4}{rgb}{0.55,0.11,0.38}
\definecolor{mediumaquamarine}{rgb}{0.40,0.80,0.67}
\definecolor{mediumblue}{rgb}{0.00,0.00,0.80}
\definecolor{mediumorchid1}{rgb}{0.88,0.40,1.00}
\definecolor{mediumorchid2}{rgb}{0.82,0.37,0.93}
\definecolor{mediumorchid3}{rgb}{0.71,0.32,0.80}
\definecolor{mediumorchid4}{rgb}{0.48,0.22,0.55}
\definecolor{mediumorchid}{rgb}{0.73,0.33,0.83}
\definecolor{mediumpurple}{rgb}{0.58,0.44,0.86}
\definecolor{mediumpurple1}{rgb}{0.67,0.51,1.00}
\definecolor{mediumpurple2}{rgb}{0.62,0.47,0.93}
\definecolor{mediumpurple3}{rgb}{0.54,0.41,0.80}
\definecolor{mediumpurple4}{rgb}{0.36,0.28,0.55}
\definecolor{medium}{rgb}{0.45,0.45,0.45}
\definecolor{mediumseagreen}{rgb}{0.24,0.70,0.44}
\definecolor{mediumsea}{rgb}{0.24,0.70,0.44}
\definecolor{mediumslateblue}{rgb}{0.48,0.41,0.93}
\definecolor{mediumslate}{rgb}{0.48,0.41,0.93}
\definecolor{mediumspringgreen}{rgb}{0.00,0.98,0.60}
\definecolor{mediumspring}{rgb}{0.00,0.98,0.60}
\definecolor{mediumturquoise}{rgb}{0.28,0.82,0.80}
\definecolor{mediumvioletred}{rgb}{0.78,0.08,0.52}
\definecolor{mediumviolet}{rgb}{0.78,0.08,0.52}
\definecolor{midnightblue}{rgb}{0.10,0.10,0.44}
\definecolor{mintcream}{rgb}{0.96,1.00,0.98}
\definecolor{mistyrose}{rgb}{1.00,0.89,0.88}
\definecolor{mistyrose1}{rgb}{1.00,0.89,0.88}
\definecolor{mistyrose2}{rgb}{0.93,0.84,0.82}
\definecolor{mistyrose3}{rgb}{0.80,0.72,0.71}
\definecolor{mistyrose4}{rgb}{0.55,0.49,0.48}
\definecolor{moccasin}{rgb}{1.00,0.89,0.71}
\definecolor{navajowhite}{rgb}{1.00,0.87,0.68}
\definecolor{navajowhite1}{rgb}{1.00,0.87,0.68}
\definecolor{navajowhite2}{rgb}{0.93,0.81,0.63}
\definecolor{navajowhite3}{rgb}{0.80,0.70,0.55}
\definecolor{navajowhite4}{rgb}{0.55,0.47,0.37}
\definecolor{navyblue}{rgb}{0.00,0.00,0.50}
\definecolor{navy}{rgb}{0.00,0.00,0.50}
\definecolor{oldlace}{rgb}{0.99,0.96,0.90}
\definecolor{olivedrab}{rgb}{0.42,0.56,0.14}
\definecolor{olivedrab1}{rgb}{0.75,1.00,0.24}
\definecolor{olivedrab2}{rgb}{0.70,0.93,0.23}
\definecolor{olivedrab3}{rgb}{0.60,0.80,0.20}
\definecolor{olivedrab4}{rgb}{0.41,0.55,0.13}
\definecolor{orange1}{rgb}{1.00,0.65,0.00}
\definecolor{orange2}{rgb}{0.93,0.60,0.00}
\definecolor{orange3}{rgb}{0.80,0.52,0.00}
\definecolor{orange4}{rgb}{0.55,0.35,0.00}
\definecolor{orangered}{rgb}{1.00,0.27,0.00}
\definecolor{orangered1}{rgb}{1.00,0.27,0.00}
\definecolor{orangered2}{rgb}{0.93,0.25,0.00}
\definecolor{orangered3}{rgb}{0.80,0.22,0.00}
\definecolor{orangered4}{rgb}{0.55,0.15,0.00}
\definecolor{orange}{rgb}{1.00,0.65,0.00}
\definecolor{orchid}{rgb}{0.85,0.44,0.84}
\definecolor{orchid1}{rgb}{1.00,0.51,0.98}
\definecolor{orchid2}{rgb}{0.93,0.48,0.91}
\definecolor{orchid3}{rgb}{0.80,0.41,0.79}
\definecolor{orchid4}{rgb}{0.55,0.28,0.54}
\definecolor{palegoldenrod}{rgb}{0.93,0.91,0.67}
\definecolor{palegreen}{rgb}{0.60,0.98,0.60}
\definecolor{palegreen1}{rgb}{0.60,1.00,0.60}
\definecolor{palegreen2}{rgb}{0.56,0.93,0.56}
\definecolor{palegreen3}{rgb}{0.49,0.80,0.49}
\definecolor{palegreen4}{rgb}{0.33,0.55,0.33}
\definecolor{paleturquoise1}{rgb}{0.73,1.00,1.00}
\definecolor{paleturquoise2}{rgb}{0.68,0.93,0.93}
\definecolor{paleturquoise3}{rgb}{0.59,0.80,0.80}
\definecolor{paleturquoise4}{rgb}{0.40,0.55,0.55}
\definecolor{paleturquoise}{rgb}{0.69,0.93,0.93}
\definecolor{palevioletred}{rgb}{0.86,0.44,0.58}
\definecolor{palevioletred1}{rgb}{1.00,0.51,0.67}
\definecolor{palevioletred2}{rgb}{0.93,0.47,0.62}
\definecolor{palevioletred3}{rgb}{0.80,0.41,0.54}
\definecolor{palevioletred4}{rgb}{0.55,0.28,0.36}
\definecolor{paleviolet}{rgb}{0.86,0.44,0.58}
\definecolor{papayawhip}{rgb}{1.00,0.94,0.84}
\definecolor{peachPuff1}{rgb}{1.00,0.85,0.73}
\definecolor{peachPuff2}{rgb}{0.93,0.80,0.68}
\definecolor{peachPuff3}{rgb}{0.80,0.69,0.58}
\definecolor{peachPuff4}{rgb}{0.55,0.47,0.40}
\definecolor{peachpuff}{rgb}{1.00,0.85,0.73}
\definecolor{peru}{rgb}{0.80,0.52,0.25}
\definecolor{pink}{rgb}{1.00,0.75,0.80}
\definecolor{pink1}{rgb}{1.00,0.71,0.77}
\definecolor{pink2}{rgb}{0.93,0.66,0.72}
\definecolor{pink3}{rgb}{0.80,0.57,0.62}
\definecolor{pink4}{rgb}{0.55,0.39,0.42}
\definecolor{plum}{rgb}{0.87,0.63,0.87}
\definecolor{plum1}{rgb}{1.00,0.73,1.00}
\definecolor{plum2}{rgb}{0.93,0.68,0.93}
\definecolor{plum3}{rgb}{0.80,0.59,0.80}
\definecolor{plum4}{rgb}{0.55,0.40,0.55}
\definecolor{powderblue}{rgb}{0.69,0.88,0.90}
\definecolor{purple}{rgb}{0.63,0.13,0.94}
\definecolor{purple1}{rgb}{0.61,0.19,1.00}
\definecolor{purple2}{rgb}{0.57,0.17,0.93}
\definecolor{purple3}{rgb}{0.49,0.15,0.80}
\definecolor{purple4}{rgb}{0.33,0.10,0.55}
\definecolor{red}{rgb}{1.00,0.00,0.00}
\definecolor{red1}{rgb}{1.00,0.00,0.00}
\definecolor{red2}{rgb}{0.93,0.00,0.00}
\definecolor{red3}{rgb}{0.80,0.00,0.00}
\definecolor{red4}{rgb}{0.55,0.00,0.00}
\definecolor{rosybrown}{rgb}{0.74,0.56,0.56}
\definecolor{rosybrown1}{rgb}{1.00,0.76,0.76}
\definecolor{rosybrown2}{rgb}{0.93,0.71,0.71}
\definecolor{rosybrown3}{rgb}{0.80,0.61,0.61}
\definecolor{rosybrown4}{rgb}{0.55,0.41,0.41}
\definecolor{royalblue}{rgb}{0.25,0.41,0.88}
\definecolor{royalblue1}{rgb}{0.28,0.46,1.00}
\definecolor{royalblue2}{rgb}{0.26,0.43,0.93}
\definecolor{royalblue3}{rgb}{0.23,0.37,0.80}
\definecolor{royalblue4}{rgb}{0.15,0.25,0.55}
\definecolor{saddlebrown}{rgb}{0.55,0.27,0.07}
\definecolor{salmon}{rgb}{0.98,0.50,0.45}
\definecolor{salmon1}{rgb}{1.00,0.55,0.41}
\definecolor{salmon2}{rgb}{0.93,0.51,0.38}
\definecolor{salmon3}{rgb}{0.80,0.44,0.33}
\definecolor{salmon4}{rgb}{0.55,0.30,0.22}
\definecolor{sandybrown}{rgb}{0.96,0.64,0.38}
\definecolor{seagreen}{rgb}{0.18,0.55,0.34}
\definecolor{seagreen1}{rgb}{0.33,1.00,0.62}
\definecolor{seagreen2}{rgb}{0.31,0.93,0.58}
\definecolor{seagreen3}{rgb}{0.26,0.80,0.50}
\definecolor{seagreen4}{rgb}{0.18,0.55,0.34}
\definecolor{seashell}{rgb}{1.00,0.96,0.93}
\definecolor{seashell1}{rgb}{1.00,0.96,0.93}
\definecolor{seashell2}{rgb}{0.93,0.90,0.87}
\definecolor{seashell3}{rgb}{0.80,0.77,0.75}
\definecolor{seashell4}{rgb}{0.55,0.53,0.51}
\definecolor{sienna}{rgb}{0.63,0.32,0.18}
\definecolor{sienna1}{rgb}{1.00,0.51,0.28}
\definecolor{sienna2}{rgb}{0.93,0.47,0.26}
\definecolor{sienna3}{rgb}{0.80,0.41,0.22}
\definecolor{sienna4}{rgb}{0.55,0.28,0.15}
\definecolor{skyblue1}{rgb}{0.53,0.81,1.00}
\definecolor{skyblue2}{rgb}{0.49,0.75,0.93}
\definecolor{skyblue3}{rgb}{0.42,0.65,0.80}
\definecolor{skyblue4}{rgb}{0.29,0.44,0.55}
\definecolor{skyblue}{rgb}{0.53,0.81,0.92}
\definecolor{slateblue1}{rgb}{0.51,0.44,1.00}
\definecolor{slateblue2}{rgb}{0.48,0.40,0.93}
\definecolor{slateblue3}{rgb}{0.41,0.35,0.80}
\definecolor{slateblue4}{rgb}{0.28,0.24,0.55}
\definecolor{slateblue}{rgb}{0.42,0.35,0.80}
\definecolor{slategray}{rgb}{0.44,0.50,0.56}
\definecolor{slategray1}{rgb}{0.78,0.89,1.00}
\definecolor{slategray2}{rgb}{0.73,0.83,0.93}
\definecolor{slategray3}{rgb}{0.62,0.71,0.80}
\definecolor{slategray4}{rgb}{0.42,0.48,0.55}
\definecolor{snow}{rgb}{1.00,0.98,0.98}
\definecolor{snow1}{rgb}{1.00,0.98,0.98}
\definecolor{snow2}{rgb}{0.93,0.91,0.91}
\definecolor{snow3}{rgb}{0.80,0.79,0.79}
\definecolor{snow4}{rgb}{0.55,0.54,0.54}
\definecolor{springgreen1}{rgb}{0.00,1.00,0.50}
\definecolor{springgreen2}{rgb}{0.00,0.93,0.46}
\definecolor{springgreen3}{rgb}{0.00,0.80,0.40}
\definecolor{springgreen4}{rgb}{0.00,0.55,0.27}
\definecolor{springgreen}{rgb}{0.00,1.00,0.50}
\definecolor{steelblue}{rgb}{0.27,0.51,0.71}
\definecolor{steelblue1}{rgb}{0.39,0.72,1.00}
\definecolor{steelblue2}{rgb}{0.36,0.67,0.93}
\definecolor{steelblue3}{rgb}{0.31,0.58,0.80}
\definecolor{steelblue4}{rgb}{0.21,0.39,0.55}
\definecolor{tan}{rgb}{0.82,0.71,0.55}
\definecolor{tan1}{rgb}{1.00,0.65,0.31}
\definecolor{tan2}{rgb}{0.93,0.60,0.29}
\definecolor{tan3}{rgb}{0.80,0.52,0.25}
\definecolor{tan4}{rgb}{0.55,0.35,0.17}
\definecolor{thistle}{rgb}{0.85,0.75,0.85}
\definecolor{thistle1}{rgb}{1.00,0.88,1.00}
\definecolor{thistle2}{rgb}{0.93,0.82,0.93}
\definecolor{thistle3}{rgb}{0.80,0.71,0.80}
\definecolor{thistle4}{rgb}{0.55,0.48,0.55}
\definecolor{tomato}{rgb}{1.00,0.39,0.28}
\definecolor{tomato1}{rgb}{1.00,0.39,0.28}
\definecolor{tomato2}{rgb}{0.93,0.36,0.26}
\definecolor{tomato3}{rgb}{0.80,0.31,0.22}
\definecolor{tomato4}{rgb}{0.55,0.21,0.15}
\definecolor{turquoise}{rgb}{0.25,0.88,0.82}
\definecolor{turquoise1}{rgb}{0.00,0.96,1.00}
\definecolor{turquoise2}{rgb}{0.00,0.90,0.93}
\definecolor{turquoise3}{rgb}{0.00,0.77,0.80}
\definecolor{turquoise4}{rgb}{0.00,0.53,0.55}
\definecolor{violetred1}{rgb}{1.00,0.24,0.59}
\definecolor{violetred2}{rgb}{0.93,0.23,0.55}
\definecolor{violetred3}{rgb}{0.80,0.20,0.47}
\definecolor{violetred4}{rgb}{0.55,0.13,0.32}
\definecolor{violetred}{rgb}{0.82,0.13,0.56}
\definecolor{violet}{rgb}{0.93,0.51,0.93}
\definecolor{wheat}{rgb}{0.96,0.87,0.70}
\definecolor{wheat1}{rgb}{1.00,0.91,0.73}
\definecolor{wheat2}{rgb}{0.93,0.85,0.68}
\definecolor{wheat3}{rgb}{0.80,0.73,0.59}
\definecolor{wheat4}{rgb}{0.55,0.49,0.40}
\definecolor{white}{rgb}{1.00,1.00,1.00}
\definecolor{whitesmoke}{rgb}{0.96,0.96,0.96}
\definecolor{yellow}{rgb}{1.00,1.00,0.00}
\definecolor{yellow1}{rgb}{1.00,1.00,0.00}
\definecolor{yellow2}{rgb}{0.93,0.93,0.00}
\definecolor{yellow3}{rgb}{0.80,0.80,0.00}
\definecolor{yellow4}{rgb}{0.55,0.55,0.00}
\definecolor{yellowgreen}{rgb}{0.60,0.80,0.20}
\makeatletter
\newlength\@tempdim@x
\newlength\@tempdim@y
\newcommand\AtUpperLeftCorner[3]{%
\begingroup
\@tempdim@x=0cm
\@tempdim@y=\paperheight
\advance\@tempdim@x#1
\advance\@tempdim@y-#2
\put(\LenToUnit{\@tempdim@x},\LenToUnit{\@tempdim@y}){#3}%
\endgroup}
\newcommand\AtUpperRightCorner[3]{%
\begingroup
\@tempdim@x=\paperwidth
\@tempdim@y=\paperheight
\advance\@tempdim@x-#1
\advance\@tempdim@y-#2
\put(\LenToUnit{\@tempdim@x},\LenToUnit{\@tempdim@y}){#3}%
\endgroup}
\newcommand\AtLowerLeftCorner[3]{%
\begingroup
\@tempdim@x=0cm
\@tempdim@y=0cm
\advance\@tempdim@x#1
\advance\@tempdim@y#2
\put(\LenToUnit{\@tempdim@x},\LenToUnit{\@tempdim@y}){#3}%
\endgroup}
\newcommand\AtLowerRightCorner[3]{%
\begingroup
\@tempdim@x=\paperwidth
\@tempdim@y=0cm
\advance\@tempdim@x-#1
\advance\@tempdim@y#2
\put(\LenToUnit{\@tempdim@x},\LenToUnit{\@tempdim@y}){#3}%
\endgroup}
\makeatother
\newtheorem{theoreme}{Theorem}[section]

\newtheorem{definition}[theoreme]{Definition}
\newtheorem{proposition}[theoreme]{Proposition}
\newtheorem{lemme}[theoreme]{Lemma}
\newtheorem{corollaire}[theoreme]{Corollary}

\newtheorem{remarque}[theoreme]{Remark}

\newenvironment{demo}{\begin{proof}}{\end{proof}}

\addto\captionsfrenchb{}
\addto\captionsfrenchb{}

    \newlength{\myarrowsize} 
    \newlength{\myoldlinewidth}

    \pgfarrowsdeclare{myto}{myto}{
        \pgfsetdash{}{0pt} 
        \pgfsetbeveljoin 
        \pgfsetroundcap 
        \setlength{\myarrowsize}{0.5pt}
        \addtolength{\myarrowsize}{.5\pgflinewidth}
        \pgfarrowsleftextend{-4\myarrowsize-.5\pgflinewidth} 
        \pgfarrowsrightextend{.7\pgflinewidth}
    }{
        \setlength{\myarrowsize}{0.4pt} 
        \addtolength{\myarrowsize}{.3\pgflinewidth}  
        \setlength{\myoldlinewidth}{\pgflinewidth}
        \pgfsetroundjoin
        \pgfsetlinewidth{0.0001pt}
        \pgfpathmoveto{\pgfpoint{0.43\myarrowsize}{0}}
        \pgfpatharc{0}{70}{0.14\myarrowsize}
        \pgfpatharc{-80}{-169.5}{4\myarrowsize}
        \pgfpatharc{150}{189}{0.95\myarrowsize and 0.95\myarrowsize}
        \pgfpatharc{0}{-40}{15\myarrowsize}
        \pgfpathmoveto{\pgfpoint{0.43\myarrowsize}{0}}
        \pgfpatharc{0}{-70}{0.14\myarrowsize}
        \pgfpatharc{80}{169.5}{4\myarrowsize}
        \pgfpatharc{-150}{-189}{0.95\myarrowsize and 0.95\myarrowsize}
        \pgfpatharc{0}{40}{15\myarrowsize}
        \pgfpathclose
        \pgfsetstrokeopacity{0.25}
        \pgfusepathqfillstroke
    }

    \pgfarrowsdeclare{myonto}{myonto}{
        \pgfsetdash{}{0pt} 
        \pgfsetbeveljoin 
        \pgfsetroundcap 
        \setlength{\myarrowsize}{0.5pt}
        \addtolength{\myarrowsize}{.5\pgflinewidth}
        \pgfarrowsleftextend{-4\myarrowsize-.5\pgflinewidth} 
        \pgfarrowsrightextend{.7\pgflinewidth}
    }{
        \setlength{\myarrowsize}{0.4pt} 
        \addtolength{\myarrowsize}{.3\pgflinewidth}  
        \setlength{\myoldlinewidth}{\pgflinewidth}
        \pgfsetroundjoin
        \pgfsetlinewidth{0.0001pt}
        \pgfpathmoveto{\pgfpoint{0.43\myarrowsize}{0}}
        \pgfpatharc{0}{70}{0.14\myarrowsize}
        \pgfpatharc{-80}{-169.5}{4\myarrowsize}
        \pgfpatharc{150}{189}{0.95\myarrowsize and 0.95\myarrowsize}
        \pgfpatharc{0}{-40}{15\myarrowsize}
        \pgfpathmoveto{\pgfpoint{-7\myarrowsize}{0}}
				\pgfpatharc{0}{70}{0.14\myarrowsize}
        \pgfpatharc{-80}{-169.5}{4\myarrowsize}
        \pgfpatharc{150}{189}{0.95\myarrowsize and 0.95\myarrowsize}
        \pgfpatharc{0}{-40}{15\myarrowsize}
        \pgfpathmoveto{\pgfpoint{0.43\myarrowsize}{0}}
        \pgfpatharc{0}{-70}{0.14\myarrowsize}
        \pgfpatharc{80}{169.5}{4\myarrowsize}
        \pgfpatharc{-150}{-189}{0.95\myarrowsize and 0.95\myarrowsize}
        \pgfpatharc{0}{40}{15\myarrowsize}
			  \pgfpathmoveto{\pgfpoint{-7\myarrowsize}{0}}
        \pgfpatharc{0}{-70}{0.14\myarrowsize}
        \pgfpatharc{80}{169.5}{4\myarrowsize}
        \pgfpatharc{-150}{-189}{0.95\myarrowsize and 0.95\myarrowsize}
        \pgfpatharc{0}{40}{15\myarrowsize}
        \pgfpathclose
        \pgfsetstrokeopacity{0.25}
        \pgfusepathqfillstroke
        \pgfusepathqfillstroke
    }

    \pgfarrowsdeclare{myhook}{myhook}{
        \setlength{\myarrowsize}{0.6pt}
        \addtolength{\myarrowsize}{.5\pgflinewidth}
        \pgfarrowsleftextend{-4\myarrowsize-.5\pgflinewidth} 
        \pgfarrowsrightextend{.7\pgflinewidth}
    }{
        \setlength{\myarrowsize}{0.6pt} 
        \addtolength{\myarrowsize}{.5\pgflinewidth}  
        \pgfsetdash{}{+0pt}
        \pgfsetroundcap
        \pgfpathmoveto{\pgfqpoint{-2pt}{-6\pgflinewidth}}
        \pgfpathcurveto
            {\pgfqpoint{4\pgflinewidth}{-4.667\pgflinewidth}}
            {\pgfqpoint{4\pgflinewidth}{0pt}}
            {\pgfpointorigin}
        \pgfusepathqstroke
    }

		\pgfarrowsdeclare{my to}{my to}
{
  \pgfarrowsleftextend{-2\pgflinewidth}
  \pgfarrowsrightextend{\pgflinewidth}
}
{
  \pgfsetlinewidth{0.8\pgflinewidth}
  \pgfsetdash{}{0pt}
  \pgfsetroundcap
  \pgfsetroundjoin
  \pgfpathmoveto{\pgfpoint{-5.5\pgflinewidth}{7.5\pgflinewidth}}
  \pgfpathcurveto
  {\pgfpoint{-4.0\pgflinewidth}{0.1\pgflinewidth}}
  {\pgfpoint{0pt}{0.25\pgflinewidth}}
  {\pgfpoint{0.75\pgflinewidth}{0pt}}
  \pgfpathcurveto
  {\pgfpoint{0pt}{-0.25\pgflinewidth}}
  {\pgfpoint{-4.0\pgflinewidth}{-0.1\pgflinewidth}}
  {\pgfpoint{-5.5\pgflinewidth}{-7.5\pgflinewidth}}
  \pgfusepathqstroke
}

\tikzstyle{vecArrow} = [thick, decoration={markings,mark=at position
   1 with {\arrow[semithick]{open triangle 60}}},
   double distance=1.4pt, shorten >= 5.5pt,
   preaction = {decorate},
   postaction = {draw,line width=1.4pt, white,shorten >= 4.5pt}]
\tikzstyle{innerWhite} = [semithick, white,line width=1.4pt, shorten >= 4.5pt]

	\makeatletter
	\newcommand\POSITION[3]{%
	\begingroup
	\@tempdim@x=0cm
	\@tempdim@y=\paperheight
	\advance\@tempdim@x#1
	\advance\@tempdim@y-#2
	\put(\LenToUnit{\@tempdim@x},\LenToUnit{\@tempdim@y}){#3}%
	\endgroup
	}
\makeatother

\addto\captionsenglish{}

\begin{document}
	\maketitle
	
	\begin{abstract}
		We define a tower  of affine Temperley-Lieb algebras of type $\tilde{A_{n}}$ and we define Markov elements in those algebras. We prove that any 	trace over an affine  Temperley-Lieb algebras of type $\tilde{A_{2}}$    is uniquely defined by its values on the Markov elements. \end{abstract}
	 
	\section{Introduction}
	
	In \cite{Sadek_Thesis} we define a tower $( \widehat{TL}_{n+1}(q))_{n\ge 0} $ of affine Temperley-Lieb algebras of type $\tilde{A_{n}}$ and we prove that there exists a unique Markov trace on this tower. 
Crucial in the proof 	is the definition of {\it Markov elements} and the following Theorem : 

	\begin{theoreme}  
	Any trace over $\widehat{TL}_{n+1} (q)$ for $2 \leq n$  is uniquely defined by its values on the Markov elements in $\widehat{TL}_{n+1} (q)$.\\		
             \end{theoreme}

The proof of this Theorem for $3 \leq n$ is given in  \cite{Sadek_2013_1}, where we have omitted  the case   $n=2$,  long and technical. We thus present it here for completeness.

		\section{Notations}
		Let $K$ be an integral domain of characteristic $0$. Suppose that $q$ is a square invertible element in $K$ of which we fix a root  $\sqrt{q}$. For $x,y$ in a given ring we define $V(x,y):= xyx+xy+yx+x+y+1$. We mean by algebra in what follows $K$-algebra.\\

		We denote by $B(\tilde{A_{n}})$ (resp. $W(\tilde{A_{n}})$) the affine braid (resp. affine Coxeter) group with $n+1$ generators of type $\tilde{A}$, while we denote by $B(A_{n})$ (resp. $W(A_{n})$) the braid (resp. Coxeter) group with $n$ generators of type $A$, where $n \geq 0 $. Let $W^{c}(\tilde{A_{n}})$ (resp. $W^{c}(A_{n}))$ be the set of fully commutative elements in $W(\tilde{A_{n}})$ (resp. $W(A_{n}))$.\\
		
		Let $ n\geq 2$. We define $\widehat{TL}_{n+1} (q)$ to be the algebra with unit given by a set of generators  $\left\{g_{\sigma_{1}}, ...,~ g_{\sigma_{n}}, g_{a_{n+1}}\right\}$, with the following relations \cite{Graham_Lehrer_1998}:\\
		
		\begin{itemize}[label=$\bullet$, font=\normalsize, font=\color{black}, leftmargin=2cm, parsep=0cm, itemsep=0.25cm, topsep=0cm] 
			\item $g_{\sigma_{i}} g_{\sigma_{j}} =g_{\sigma_{j}} g_{\sigma_{i}} $,~~~~~~~~~~~~~~~~~~~~~~~~~~~~for  $1\leq i,j\leq n$ and $ \left| i-j\right| \geq 2$. 
			\item $g_{\sigma_{i}} g_{a_{n+1}} =g_{a_{n+1}} g_{\sigma_{i}} $,~~~~~~~~~~~~~~~~~~~~~~~for  $2\leq i \leq n-1$.
			\item $g_{\sigma_{i}}g_{\sigma_{i+1}}g_{\sigma_{i}} = g_{\sigma_{i+1}}g_{\sigma_{i}}g_{\sigma_{i+1}}$,~~~~~~~~~~~~~~for $1\leq i\leq n-1$.
			\item $g_{\sigma_{i}}g_{a_{n+1}}g_{\sigma_{i}} = g_{a_{n+1}}g_{\sigma_{i}}g_{a_{n+1}}$,~~~~~~~~~~~~~for $i= 1, n $.
			\item $g^{2}_{\sigma_{i}} = (q-1)g_{\sigma_{i}} +q$, ~~~~~~~~~~~~~~~~~~~~~for $1\leq i\leq n$.
			\item $g^{2}_{a_{n+1}} = (q-1)g_{a_{n+1}} +q$, ~~~~~~~~~
			\item $V(g_{\sigma_{i}},g_{\sigma_{i+1}})=V(g_{\sigma_{1}},g_{a_{n+1}}) = V(g_{\sigma_{n}},g_{a_{n+1}})= 0$, ~~~~~~~~~~ for $1\leq i\leq n-1$.\\
		\end{itemize}
		
		 The set $\left\{ g_{w}: w \in W^{c}(\tilde{A_{n}})\right\}$ is well defined in the usual sense of the theory of Hecke algebra and it is a $K$-basis. We set $T_{a_{n+1}}$ (resp. $T_{\sigma_{i}}$ for $1 \leq i \leq n$) to be $\sqrt{q}g_{a_{n+1}}$ (resp. $\sqrt{q}g_{\sigma_{i}}$ for $1 \leq i \leq n$). Hence, $T_{w}$ is well defined for $w \in W^{c}(\tilde{A_{n}})$, it equals $q^{\frac{l(w)}{2}}g_{w}$. The multiplication associated to the basis  $\left\{ T_{w}: w \in W^{c}(\tilde{A_{n}})\right\}$, is given as follows:
		 
		\begin{eqnarray}
			T_{w} T_{v} &=& T_{wv} ~~~~~~~~~~~~~~~~~~~~~~~~~~~~~~ \text{whenever } l(wv) = l(w) + l (v).\nonumber\\
			T_{s} T_{w} &=& \sqrt{q}(q-1) T_{w} +q^{2}  T_{sw} ~~~~~~~~ \text{whenever } l(sw) = l(w) - 1, \nonumber
		\end{eqnarray}

		for $w,v$ in $W^{c}(\tilde{A_{n}})$ and $s$ in $\left\{\sigma_{1},...,~\sigma_{n}, a_{n+1}\right\}$.\\
		
		 In what follows we suppose that $q+1$ is  invertible in $K$, we set $\delta = \frac{1}{2+q+q^{-1}} = \frac{q}{(1+q)^{2}}$ in $K$. In view of \cite{Graham_Lehrer_2003} , for $1 \leq i \leq n$ we set $ f_{\sigma_{i}}:= \frac{ g_{\sigma_{i}}+1}{q+1}$ and $ f_{a_{n+1}}:= \frac{g_{a_{n+1}}+1}{q+1}$. In other terms $ g_{\sigma_{i}} = (q+1) f_{\sigma_{i}} -1$, and $ g_{a_{n+1}} = (q+1) f_{a_{n+1}} -1$. The set $\left\{ f_{w}: w \in W^{c}(\tilde{A_{n}})\right\}$ is well defined and it is a $K$-basis for $\widehat{TL}_{n+1} (q)$.\\
		
		We define the Temperley-Lieb algebra of type $A$ with $n$ generators $TL_{n}(q)$, as the subalgebra of $\widehat{TL}_{n+1} (q)$ generated by  $\left\{g_{\sigma_{1}} ,...,~ g_{\sigma_{n}}\right\}$, with $\left\{ g_{w}: w \in W^{c}(A_{n})\right\}$ as $K$-basis.\\
		
		Now for $TL_{0}(q) = K$, we consider the following tower: 
		
		\begin{eqnarray}
				TL_{0}(q)~\subset TL_{1}(q) ~~...\subset TL_{n-1}(q) ~~\subset TL_{n}(q) ~~... \nonumber\\\nonumber
			\end{eqnarray}
			
		\begin{theoreme} \label{1_1}\cite{Jones_1985}	
		There is a unique collection of traces $(\tau_{n+1})_{0 \leq n}$ on $(TL_{n})_{0\leq n}$, such that: \\
			\begin{itemize}[label=$\bullet$, font=\normalsize, font=\color{black}, leftmargin=2cm, parsep=0cm, itemsep=0.25cm, topsep=0cm]
				\item[$(1)$] $\tau_{1}(1) = 1 $.
				\item[$(2)$] For $ 1 \leq n$, we have $\tau_{n+1}(hT^{\pm1}_{\sigma_{n}}) = \tau_{n}(h) $, for any $h$ in $TL_{n-1}(q)$.
			\end{itemize}	
		\end{theoreme}		
		
					\vspace{0.25cm}
		The collection  $(\tau_{n+1})_{0 \leq n}$ is called a Markov trace. Moreover, for any $a,b$ and $c$ in $TL_{n}(q)$ and for $ n \geq 1$, every $\tau_{n+1}: TL_{n}(q) \longrightarrow K$ verifies: 
					\begin{eqnarray}
						\tau_{n+1}(bT_{\sigma_{n}}c)= \tau_{n}(bc)$ and $\tau_{n+1}(a)=- \frac{1+q}{\sqrt{q}}\tau_{n}(a). \nonumber
					\end{eqnarray}
			
	 	\clearpage

	\section{The tower of affine Temperley-Lieb algebras and affine Markov trace}
	
	In this section we define a tower of affine Temperley-Lieb algebras, we show that this tower "surjects" onto the tower of Temperley-Lieb algebras mentioned in the introduction, and we define the affine Markov trace.\\ 
	
	We consider the Dynkin diagram of the group $B(\tilde{A_{n}})$. We denote the Dynkin automorphism $(\sigma_{1}\mapsto \sigma_{2} \mapsto ... \sigma_{n} \mapsto a_{n+1} \mapsto \sigma_{1})$ by $\psi_{n+1}$. Notice that $ \sigma_{n} \sigma_{n-1}.. \sigma_{1}a_{n+1} $ acts on $ B(\tilde{A_{n-1}}) $ as $\psi_{n}$ as follows ($\sigma_{1} \mapsto a_{n} \mapsto \sigma_{n-1} \mapsto \sigma_{n-2} \mapsto.. \sigma_{2} \mapsto \sigma_{1}$). We write $(\sigma_{n} .. \sigma_{1} a_{n+1})^{d}h = \psi^{d} \left[ h\right]  (\sigma_{n} .. \sigma_{1} a_{n+1})^{d} $, for any $h$ in $ B(\tilde{A_{n-1}}) $, we keep same convention for the affine Temperley-Lieb algebra.	\\
			
					
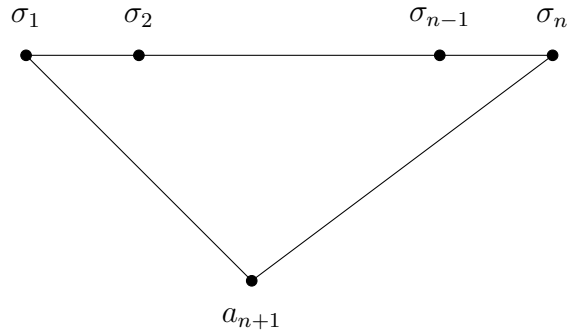
			\begin{figure}[ht]
				\centering
				\begin{tikzpicture}

 \node at (0,0.5) {$\sigma_{1}$}; 
  \filldraw (0,0) circle (2pt);
   
  \draw (0,0) -- (1.5, 0);
  
  \node at (1.5,0.5) {$\sigma_{2}$};
  \filldraw (1.5,0) circle (2pt);

  \draw (1.5,0) -- (5.5, 0);

  \node at (5.5,0.5) {$\sigma_{n-1}$};
  \filldraw (5.5,0) circle (2pt);
 
  \draw (5.5,0) -- (7, 0);
  
  \node at (7,0.5) {$\sigma_{n}$};
  \filldraw (7,0) circle (2pt);

  \draw (7,0) -- (3, -3);
  
  \filldraw (3, -3) circle (2pt);
  \node at (3, -3.5) {$a_{n+1}$};

  \draw (3, -3) -- (0, 0);
               \end{tikzpicture}
			\end{figure}
	
	We have the following injection
	
	\begin{eqnarray}
					G_{n}: K[B(\tilde{A_{n-1}})] &\longrightarrow& K[B(\tilde{A_{n}})]  \nonumber\\
					\sigma_{i} &\longmapsto& \sigma_{i}$ ~~~ \text{for} $1\leq i\leq n-1 \nonumber\\
					a_{n} &\longmapsto& \sigma_{n} a_{n+1}\sigma^{-1}_{n} \nonumber
				\end{eqnarray}
	
	We prove in \cite{Sadek_Thesis}, to which we refer for details, the following two propositions:\\
	
	
	\begin{proposition} \label{2_1} The injection $G_{n}$ induces the following morphism of algebras:
	
			\begin{eqnarray}
				F_{n}: \widehat{TL}_{n}(q) &\longrightarrow& \widehat{TL}_{n+1}(q) \nonumber\\
				t_{\sigma_{i}} &\longmapsto & g_{\sigma_{i}} \text{ for } 1 \leq i \leq n-1 \nonumber\\
				t_{a_{n}} &\longmapsto & g_{\sigma_{n}} g_{a_{n+1}} g^{-1}_{\sigma_{n}}. \nonumber
			\end{eqnarray}
	\end{proposition}	
		
	\begin{proposition} The following map is a surjection of algebras
	
			\begin{eqnarray}
				E_{n}: \widehat{TL}_{n+1}(q) &\longrightarrow& TL_{n}(q) \nonumber\\
				g_{\sigma_{i}} &\longmapsto & g_{\sigma_{i}} \text{ for } 1 \leq i \leq n \nonumber\\
				g_{a_{n+1}} &\longmapsto & g_{\sigma_{1}} ... g_{\sigma_{n-1}} g_{\sigma_{n}} g^{-1}_{\sigma_{n-1}} ... g^{-1}_{\sigma_{1}}. \nonumber
			\end{eqnarray}
	
	Moreover, the following diagram commutes: 
	
	\begin{center}    

	\begin{tikzpicture}

			\matrix[matrix of math nodes,row sep=1cm,column sep=1cm]{
			|(A)| \widehat{TL}_{n} (q)    & & & &    |(B)| \widehat{TL}_{n+1} (q)  \\
			                              & & & &                      \\								
			|(C)| TL_{n-1}(q)             & & & &    |(D)| TL_{n}(q)    \\
				};
				
				\path (C) edge[-myhook,line width=0.42pt] node[above, xshift=-5mm, yshift=-2mm, rotate=0] {\footnotesize $E_{n-1}$}   (A);
				\path (A) edge[-myonto,line width=0.42pt]  node[above, xshift=-5mm, yshift=-2mm, rotate=0] {\footnotesize $ $}    (C);
				
				\path (A) edge[-myto,line width=0.42pt] node[above, xshift=1.5mm, yshift=0mm, rotate=0] {\footnotesize $F_{n}$}      (B);
				
				\path (D) edge[-myhook,line width=0.42pt] node[below, xshift=1.5mm, yshift=0mm, rotate=0] {\footnotesize $$}    (C);
				\path (C) edge[-myto,line width=0.42pt]      (D);
								
				\path (B) edge[-myonto,line width=0.42pt]   node[above, xshift=5mm, , yshift=-2mm, rotate=0] {\footnotesize $E_{n} $}   (D);

		\end{tikzpicture}	
	\end{center}
	
	 Moreover, it is immediate that $E_{n}$ composed with the natural inclusion of $TL_{n}(q)$ into $\widehat{TL}_{n+1}(q)$, gives $Id_{TL_{n}(q)}$.
	\end{proposition}

			In view of proposition \ref{2_1} we can consider the tower of affine T-L algebras (it is not known whether it is a tower of faithful arrows or not):\\ 
		\begin{eqnarray}
			\widehat{TL}_{1}(q) \stackrel{F_{1}}{\longrightarrow}  \widehat{TL}_{2}(q) \stackrel{F_{2}} {\longrightarrow}\widehat{TL}_{3}(q) \longrightarrow ~~...~~ \widehat{TL}_{n}(q) \stackrel{F_{n}} {\longrightarrow}\widehat{TL}_{n+1}(q) \longrightarrow ... \nonumber\\\nonumber
		\end{eqnarray}
				
			\begin{definition} \label{5_2_1}
				We call $(\hat{\tau}_{n})_{1 \leq n}$ an affine Markov trace, if every $\hat{\tau}_{n}$ is a trace function on $\widehat{TL}_{n} (q)$ with the following conditions:\\
			
				\begin{itemize}
					\item $\hat{\tau}_{1}(1) = 1$, (here $\widehat{TL}_{1} (q) = K$).				
					\item $\hat{\tau}_{n+1}(F_{n}(h)T^{\pm1}_{\sigma_{n}}) =  \hat{\tau}_{n}(h)$, for all $h \in  \widehat{TL}_{n} (q)$ and for $n \geq 1$.
					\item $\hat{\tau}_{n}$ is invariant under the Dynkin automorphism $\psi_{n}$ for all $n$.\\               
				\end{itemize}
			\end{definition}
		
		\begin{remarque} 
				We notice that the second condition gives us that $\hat{\tau}_{n+1}\big(F_{n}(h)T^{-1}_{\sigma_{n}}\big) =  \hat{\tau}_{n}\big(h\big)$, which means that:
				\begin{eqnarray}
					\hat{\tau}_{n+1}\big(F_{n}(h)[\frac{1}{q^{2}}T_{\sigma_{n}}- \frac{q-1}{q\sqrt{q}}] \big) =  \hat{\tau}_{n}\big(h). \text{ Thus } \hat{\tau}_{n+1} \big(F_{n}(h)\big) = -\frac{q+1}{\sqrt{q}} \hat{\tau}_{n}\big(h\big). \nonumber\\\nonumber
				\end{eqnarray}
			\end{remarque}
			
			\begin{remarque}
			
			The third condition of definition \ref{5_2_1} is, in fact, not independent, i.e., it results from the first and second conditions (see \cite{Sadek_2013_1}).   Nevertheless, we will keep viewing it as a condition. \\
	
	      \end{remarque}
		
		\begin{remarque}
		  This affine Markov trace does the job topologically, i.e., it gives an invariant for "affine oriented knots" and generalizes, in fact, the Jones invariant, noticing that the set of oriented knots in $S^{3}$ injects naturally into the set of "affine oriented knots". For further details see \cite{Sadek_Thesis}.

		\end{remarque}

			Now, consider the following commutative diagram: \\
			
			\begin{tikzpicture}

			\matrix[matrix of math nodes,row sep=1.3cm,column sep=0.5cm]{
			|(A)|\widehat{TL}_{1}(q) & & & & |(B)| \widehat{TL}_{2}(q) & & |(BB)| \textcolor{white}{T}\dots ~~~~ & &|(C)|\widehat{TL}_{n}(q) & & & & |(D)| \widehat{TL}_{n+1}(q) \\
			|(E)|          TL_{0}(q) & & & & |(F)|           TL_{1}(q) & & |(FF)| \textcolor{white}{T}\dots ~~~~ & &|(G)| TL_{n-1}(q)        & & & & |(H)|             TL_{n}(q) \\
				                    & & & &                           & &                                       & &                         & & & &                             \\
				                    & & & &                           & &|(I) |K                                & &                         & & & &                             \\
				};

			\path (A) edge[-myto,line width=0.42pt]                                                                             (B);

			\path (B) edge[-myto,line width=0.42pt]                                                                             (BB);

			\path (BB) edge[-myto,line width=0.42pt]                                                                             (C);
			
			\path (C) edge[-myto,line width=0.42pt]                                                                             (D);
			
			
			\path (F) edge[-myhook,line width=0.42pt]                                                                           (E);
			\path (E) edge[-myto,line width=0.42pt]                                                                             (F);

			\path (FF) edge[-myhook,line width=0.42pt]                                                                           (F);
			\path (F) edge[-myto,line width=0.42pt]                                                                             (FF);
			
			\path (G) edge[-myhook,line width=0.42pt]                                                                           (FF);
			\path (FF) edge[-myto,line width=0.42pt]                                                                             (G);

			\path (H) edge[-myhook,line width=0.42pt]                                                                           (G);
			\path (G) edge[-myto,line width=0.42pt]                                                                             (H);


			\path (A) edge[-myonto,line width=0.42pt]                                                                           (E);
			
			\path (B) edge[-myonto,line width=0.42pt]                                                                           (F);

			\path (C) edge[-myonto,line width=0.42pt]                                                                           (G);

			\path (D) edge[-myonto,line width=0.42pt]                                                                           (H);

			\path (E) edge[-myto,line width=0.42pt]  node[above, xshift=-5mm, yshift=-2mm, rotate=0] {\footnotesize $\tau_{1}$}  (I);

			\path (F) edge[-myto,line width=0.42pt]  node[above, xshift=-3mm, yshift=-2mm, rotate=0] {\footnotesize $\tau_{2}$}  (I);

			\path (G) edge[-myto,line width=0.42pt]  node[above, xshift=3mm, yshift=-2mm, rotate=0] {\footnotesize $\tau_{n}$}  (I);

			\path (H) edge[-myto,line width=0.42pt]  node[above, xshift=8mm, yshift=-2mm, rotate=0] {\footnotesize $\tau_{n+1}$}  (I);

		\end{tikzpicture}

			Set $\rho_{n+1}$ to be the trace over  $\widehat{TL}_{n+1} (q) $ induced by $\tau_{n+1}$ over  $TL_{n} (q) $ for $0 \leq n$.  	
		 We prove in   		 \cite{Sadek_2013_1} that  $(\rho_{i})_{1 \leq i}$ is an affine Markov trace over $\big( \widehat{TL}_{i} (q)\big)_{1 \leq i}$ and we prove   the following Theorem:  
		
			\begin{theoreme}\label{5_2_8} \cite{Sadek_2013_1} 

		There exists a unique affine Markov trace over the tower of $\tilde{A}$-type Temperley-Lieb algebras, namely $(\rho_{i})_{1 \leq i}$.   
	\end{theoreme}
		
	The proof relies on Theorem 	 \ref{5_1_1}  below, the proof of which separates into two  cases: $n=2$ and $n \ge 3$. The latter case is included in 
	\cite{Sadek_2013_1} while the former  appears in the present note. 
		
	
	\section{Markov elements and traces on $\widehat{TL}_{n+1}(q)$}

	 			\subsection{Markov elements}
		
		We consider $F_{n}:\widehat{TL}_{n} (q)\longrightarrow \widehat{TL}_{n+1} (q)$ of proposition \ref{2_1}. In this subsection we set $F:= F_{n}$. We give a definition of Markov elements in $\widehat{TL}_{n+1} (q)$ for $ 2 \leq n $. Then we show that any trace over $\widehat{TL}_{n+1} (q)$ is uniquely determined by its values on those elements.

			\begin{definition}
				For $F$ as above, and $n \geq 2$, a Markov element in $\widehat{TL}_{n+1} (q)$ is any element of the form $ A g^{\epsilon}_{\sigma_{n}} B $, where $A$ and $B$ are in $F(\widehat{TL}_{n} (q))$ and $\epsilon \in \left\{ 0,1 \right\}$. 
			\end{definition}
		 
		The aim of this subsection is to prove the following theorem for $n=2$.
			
			\vspace{0.25cm}		
			
			\begin{theoreme}	\cite{Sadek_2013_1}  \label{5_1_1} 
				Let $\tau_{n+1}$ be any trace over $\widehat{TL}_{n+1} (q)$ for $2 \leq n$. Then, $\tau_{n+1}$ is uniquely defined by its values on the Markov elements in $\widehat{TL}_{n+1} (q)$.\\		
             \end{theoreme}

			\vspace{0.25cm}
			The proof of theorem \ref{5_1_1} for $n=2$ is divided into two parts. 
			In the first we show some general facts, in the second we give the proof   for $ n=2$.\\	
		
			\clearpage
			
			\textbf{Part 1 }\\
	
	In this part, we suppose that $\tau_{n+1}$ is any trace on $\widehat{TL}_{n+1} (q)$. We will apply $\tau_{n+1}$ to $\widehat{TL}_{n+1} (q)$ assuming that $2 \leq n$, and show that  $\tau_{n+1}$ is uniquely determined on $\widehat{TL}_{n+1} (q)$ by its values on the positive powers of $g_{\sigma_{n} \sigma_{n-1} ..\sigma_{1} a_{n+1}}$, in addition to its values on Markov elements. From now on we denote by $w$: an arbitrary element in $W^{c}(\tilde{A_{n}})$.
	
			\begin{lemme}\label{5_1_3}
				In $\widehat{TL}_{n+1} (q) $ we have: \\

					\vspace{-0.8cm}
					\begin{eqnarray}
						(1)~g_{\sigma_{n}} (g_{\sigma_{n} \sigma_{n-1} ..\sigma_{1} a_{n+1}})^{k} &=& (q-1) (g_{\sigma_{n} \sigma_{n-1} ..\sigma_{1} a_{n+1}})^{k} + \sum\limits^{i=k-1}_{i=1} f_{i}  (g_{\sigma_{n} \sigma_{n-1} ..\sigma_{1} a_{n+1}})^{i} \nonumber \\
						& & + A \big(g_{\sigma_{n-1} \sigma_{n-2} ..\sigma_{1}}F(t_{a_{n}})\big)^{k} g_{\sigma_{n}}\prod^{j=k-1}_{j=0} \psi^{j} \big[F((t_{a_{n}})^{-1})\big], \nonumber
					\end{eqnarray}

					\vspace{-0.5cm}
					\begin{eqnarray}
						(2)~(g_{\sigma_{n} \sigma_{n-1} ..\sigma_{1} a_{n+1}})^{k}g_{\sigma_{n}} &=& (q-1) (g_{\sigma_{n} \sigma_{n-1} ..\sigma_{1} a_{n+1}})^{k} + \sum^{i=k-1}_{i=1} h_{i}  (g_{\sigma_{n} \sigma_{n-1} ..\sigma_{1} a_{n+1}})^{i}\nonumber \\
						& & + A \prod^{j=k-1}_{j=0} \phi^{j} \big[(g_{\sigma_{{n-1}}})^{-1}\big] g_{\sigma_{n}}\big(g_{\sigma_{n-1} \sigma_{n-2} ..\sigma_{1}}F(t_{a_{n}})\big)^{k}, \nonumber
					\end{eqnarray}
	
				with $ A $ in the ground field, $ f_{i},h_{i} $ in $F (\widehat{TL}_{n} (q)) $ and $ \phi^{-1} = \psi $.
					
			\end{lemme}
			
			\begin{demo}
				\begin{eqnarray}
					g_{\sigma_{n}} \big(g_{\sigma_{n} \sigma_{n-1} ..\sigma_{1} a_{n+1}}\big)^{k} &=& \big(q-1\big) \big(g_{\sigma_{n} \sigma_{n-1} ..\sigma_{1} a_{n+1}}\big)^{k} \nonumber\\
					& & + q g_{\sigma_{n-1} \sigma_{n-2} ..\sigma_{1}}F\big(t_{a_{n}}\big) g_{\sigma_{n}} F\big((t_{a_{n}})^{-1}\big) \big(g_{\sigma_{n} \sigma_{n-1} ..\sigma_{1} a_{n+1}}\big)^{k-1} \nonumber\\\nonumber\\
					&=& \big(q-1\big) \big(g_{\sigma_{n} \sigma_{n-1} ..\sigma_{1} a_{n+1}}\big)^{k} + \nonumber\\
					& & q g_{\sigma_{n-1} \sigma_{n-2} ..\sigma_{1}}F\big(t_{a_{n}}\big) g_{\sigma_{n}} \big(g_{\sigma_{n} \sigma_{n-1} ..\sigma_{1} a_{n+1}}\big)^{k-1}  \psi^{k-1} \big[F((t_{a_{n}})^{-1})\big]. \nonumber\\\nonumber
				\end{eqnarray}
							
				So, by induction on $k$, (1) follows. In the very same way we deal with (2), by noticing that: $g_{a_{n+1}}g_{\sigma_{n}} =g^{-1}_{\sigma_{n}} F(t_{a_{n}}) g^{2}_{\sigma_{n}} = (q-1) g_{a_{n+1}} + qg^{-1} _{\sigma_{n}} F(t_{a_{n}}) $.

			\end{demo}
	
 			A main result in \cite{Sadek_2013_2} is to give a general form for ``fully commutative braids'', from which we deduce that any element of the basis of $\widehat{TL}_{n+1} (q)$ (where we have the convention $\sigma_{n+1} = 1$ in $W(\tilde{A_{n}})$ thus $ g_{\sigma_{n} \sigma_{n-1} .. \sigma_{i}} =1$ when $i=n+1$), is either of the form\\
 			$$ c (g_{\sigma_{n} \sigma_{n-1} ..\sigma_{1} a_{n+1}})^{k} g_{\sigma_{n}\sigma_{n-1} .. \sigma_{i}} $$
 			
 			or of the form\\
 			
 			$$g_{\sigma_{i_{0}} .. \sigma_{2}\sigma_{1}a_{n+1}} (g_{\sigma_{n} \sigma_{n-1} ..\sigma_{1} a_{n+1}})^{k} d  g_{\sigma_{n} \sigma_{n-1} .. \sigma_{i}} $$\\
 			
 			where $c$ and $d$ are in $F (\widehat{TL}_{n} (q)) $, $ 1\leq i \leq n+1 $ and $ 0 \leq i_{0} \leq n-1 $ . \\
	
			By lemma \ref{5_1_3} $ c (g_{\sigma_{n} \sigma_{n-1} ..\sigma_{1} a_{n+1}})^{k} g_{\sigma_{n}\sigma_{n-1} .. \sigma_{i}} $ is of the form:
			\begin{eqnarray}
				\sum^{j=h}_{j=1}c_{j}  (g_{\sigma_{n} \sigma_{n-1} ..\sigma_{1} a_{n+1}})^{j} + M. \nonumber
			\end{eqnarray}
			
			Where $h \leq k $, $c_{j}$ is in  $F (\widehat{TL}_{n} (q)) $ for any $ j $ and $M$ is a Markov element.\\ 
			
			Now we deal with the second form:
			\begin{eqnarray}
				\tau_{n+1} &\big(&g_{\sigma_{i_{0}} .. \sigma_{2}\sigma_{1}a_{n+1}} c (g_{\sigma_{n} \sigma_{n-1} ..\sigma_{1} a_{n+1}})^{k}  g_{\sigma_{n} \sigma_{n-1} .. \sigma_{i}}\big) = \tau_{n+1} \big(g_{\sigma_{n} \sigma_{n-1} .. \sigma_{i}} g_{\sigma_{i_{0}} .. \sigma_{2}\sigma_{1}a_{n+1}} c (g_{\sigma_{n} \sigma_{n-1} ..\sigma_{1} a_{n+1}})^{k}\big). \nonumber
			\end{eqnarray}
	
			For any possible value for $ i_{0} $ or $i$, we see that:
			\begin{eqnarray}
				g_{\sigma_{n} \sigma_{n-1} .. \sigma_{i}} g_{\sigma_{i_{0}} .. \sigma_{2}\sigma_{1}a_{n+1}} c (g_{\sigma_{n} \sigma_{n-1} ..\sigma_{1} a_{n+1}})^{k} = c'g_{\sigma_{n}} (g_{\sigma_{n} \sigma_{n-1} ..\sigma_{1} a_{n+1}})^{s} c'' ,\nonumber
			\end{eqnarray}
			
			where $ c',c''$ are in $F (\widehat{TL}_{n} (q)) $ and $ s \leq k+1 $. By lemma \ref{5_1_3} we see that this element is of the form:
			\vspace{-0.25cm}  
			\begin{eqnarray}
				\sum^{j=h}_{j=1}f_{j}  (g_{\sigma_{n} \sigma_{n-1} ..\sigma_{1} a_{n+1}})^{j} + M, \nonumber
			\end{eqnarray}
			
			where $h \leq k+1 $, $f_{j}$ is in  $F (\widehat{TL}_{n} (q)) $ for any $ j $ and $M$ is a Markov element .\\ 
	        
			Hence, we see that in order to define $\tau_{n+1}$ uniquely it is enough to have its values on Markov elements and its values on $\Omega (g_{\sigma_{n} \sigma_{n-1} ..\sigma_{1} a_{n+1}})^{k}$, where $1 \leq k$ (since if $k$ is equal to 0 then we are again in the case of a Markov element) and $\Omega$ is in $F\big(\widehat{TL}_{n} (q)\big)$.
			\vspace{0.25cm}
			\begin{lemme}\label{5_1_4}
				Let $2 \leq n $ then  $\tau_{n+1}$ is uniquely defined by its values on Markov elements, in addition to its values on $ (g_{\sigma_{n} \sigma_{n-1} ..\sigma_{1} a_{n+1}} )^{k} $,  with $0 \leq k $ .\\
			\end{lemme}

			\begin{demo} 
		 		In order to determine $\tau_{n+1} \big( h(g_{\sigma_{n} \sigma_{n-1} ..\sigma_{1} a_{n+1}})^{k}\big)$, with a positive $k$ and an arbitrary $h$ in  $F\big(\widehat{TL}_{n} (q)\big)$, it is enough to treat $\tau_{n+1} \big(F(t_{x})(g_{\sigma_{n} \sigma_{n-1} ..\sigma_{1} a_{n+1}})^{k}\big)$, with $x$ in $W^{c}(\tilde{A_{n-1}})$, but the fact that $\tau_{n+1}$ is a trace, in addition to the fact that $g_{\sigma_{n} \sigma_{n-1} ..\sigma_{1} a_{n+1}}$ acts as a Dynkin automorphism on $F \big(\widehat{TL}_{n} (q)\big)$, authorizes us to suppose that $x$ has a reduced expression which ends with $\sigma_{n-1}$. \\
		 		
		 		Now we show by induction on $l(x)$, that $\tau_{n+1} \big( F(t_{x}) (g_{\sigma_{n} \sigma_{n-1} ..\sigma_{1} a_{n+1}})^{k}\big)$ is a sum of values of $\tau_{n+1}$ over $ (g_{\sigma_{n} \sigma_{n-1} ..\sigma_{1} a_{n+1}})^{k}$, elements of the form $ h (g_{\sigma_{n} \sigma_{n-1} ..\sigma_{1} a_{n+1}})^{i}$ with $i < k $ and Markov elements, (of course with coefficients in the ground ring which might be zeros).\\ 
	
				For $l(x) = 0$ the property is true. Take $l(x) > 0$, and let $ x = z \sigma_{n-1}$ be a reduced expression, hence:
				\begin{eqnarray}
					\tau_{n+1} \big( F(t_{x}) (g_{\sigma_{n} \sigma_{n-1} ..\sigma_{1} a_{n+1}})^{k}\big) &=& \tau_{n+1}\big( F(t_{z}) F(t_{\sigma_{n-1}}) g_{\sigma_{n} \sigma_{n-1} ..\sigma_{1} a_{n+1}}(g_{\sigma_{n} \sigma_{n-1} ..\sigma_{1} a_{n+1}})^{k-1}\big) \nonumber\\\nonumber\\
					&=& \tau_{n+1} \big(F(t_{z}) \underbrace{g_{\sigma_{n-1}} g_{\sigma_{n}} g_{\sigma_{n-1}}}_{=-V(g_{\sigma_{n-1}}, g_{\sigma_{n}})} g_{ \sigma_{n-2} ..\sigma_{1} a_{n+1}}(g_{\sigma_{n} \sigma_{n-1} ..\sigma_{1} a_{n+1}})^{k-1}\big). \nonumber
				\end{eqnarray}
	
				Recalling that $V(g_{\sigma_{n-1}}, g_{\sigma_{n}})=0$, this is equal to the following sum:
				\begin{eqnarray}
					& & - \tau_{n+1}\big( F(t_{z}) (g_{\sigma_{n} \sigma_{n-1} ..\sigma_{1} a_{n+1}})^{k}\big) \nonumber\\\nonumber\\
					& & - \tau_{n+1}\big( F(t_{z}) g_{\sigma_{n-1}}   g_{ \sigma_{n-2} ..\sigma_{1}}  g_{a_{n+1}}(g_{\sigma_{n} \sigma_{n-1} ..\sigma_{1} a_{n+1}})^{k-1}\big)\nonumber\\\nonumber\\
					& & - \tau_{n+1}\big( F(t_{z}) g_{ \sigma_{n-2} ..\sigma_{1} a_{n+1}}(g_{\sigma_{n} \sigma_{n-1} ..\sigma_{1} a_{n+1}})^{k-1}\big)\nonumber\\\nonumber\\
					& & - \tau_{n+1}\big( F(t_{z}) g_{\sigma_{n-1}}   g_{ \sigma_{n-2} ..\sigma_{1}}\underbrace{ g_{\sigma_{n}} g_{a_{n+1}}}_{}(g_{\sigma_{n} \sigma_{n-1} ..\sigma_{1} a_{n+1}})^{k-1}\big) \nonumber\\
					& & - \tau_{n+1}\big( F(t_{z}) g_{ \sigma_{n-2} ..\sigma_{1}} \underbrace{g_{\sigma_{n}} g_{a_{n+1}}}_{}(g_{\sigma_{n} \sigma_{n-1} ..	\sigma_{1} a_{n+1}})^{k-1}\big).\nonumber
				\end{eqnarray}
				
				Now we apply the induction hypothesis to the first term. The second and the third terms are equal to:
				\begin{eqnarray}
					\tau_{n+1} &\bigg(& F\big(t_{z}\big) g_{\sigma_{n-1}} g_{ \sigma_{n-2} ..\sigma_{1}} F\big( t_{a_{n}}\big) g_{\sigma_{n}} F\big( (t_{a_{n}})^{-1}\big)\big(g_{\sigma_{n} \sigma_{n-1} ..\sigma_{1} a_{n+1}}\big)^{k-1}\bigg)\nonumber\\\nonumber\\
					&+& \tau_{n+1} \bigg( F\big(t_{z}\big) g_{\sigma_{n-2} ..\sigma_{1}} F\big( t_{a_{n}}\big) g_{\sigma_{n}}F\big( (t_{a_{n}})^{-1}\big)\big(g_{\sigma_{n} \sigma_{n-1} ..\sigma_{1} a_{n+1}}\big)^{k-1}\bigg), \nonumber					
				\end{eqnarray}
				
				which is equal to:
				\begin{eqnarray}
					\tau_{n+1} &\bigg(& \psi^{1-k} \big[ F\big( (t_{a_{n}})^{-1}\big) \big]F\big(t_{z}\big) g_{\sigma_{n-1}}   g_{ \sigma_{n-2} ..\sigma_{1}} F\big( t_{a_{n}}\big) \big(g_{\sigma_{n}} (g_{\sigma_{n} \sigma_{n-1} ..\sigma_{1} a_{n+1}})^{k-1}\big)\bigg)\nonumber\\\nonumber\\	
					&+& \tau_{n+1}\bigg(\psi^{1-k} \big[ F\big( (t_{a_{n}})^{-1}\big) \big] F\big(t_{z}\big) g_{ \sigma_{n-2} ..\sigma_{1}} F( t_{a_{n}})\big( g_{\sigma_{n}}(g_{\sigma_{n} \sigma_{n-1} ..\sigma_{1} a_{n+1}})^{k-1}\big)\bigg).\nonumber				
				\end{eqnarray}
				
				The fourth and the fifth terms are equal to:
				\begin{eqnarray}
					\tau_{n+1} &\bigg(& F\big(t_{z}\big) g_{\sigma_{n-1}} g_{ \sigma_{n-2} ..\sigma_{1}} F\big( t_{a_{n}}\big) \big(g_{\sigma_{n}}(g_{\sigma_{n} \sigma_{n-1} ..\sigma_{1} a_{n+1}})^{k-1}\big)\bigg) \nonumber\\\nonumber\\
					&+& \tau_{n+1}\bigg( F\big(t_{z}\big) g_{\sigma_{n-2} ..\sigma_{1}} F\big(t_{a_{n}}\big) \big(g_{\sigma_{n}}(g_{\sigma_{n} \sigma_{n-1} ..\sigma_{1} a_{n+1}})^{k-1}\big)\bigg).\nonumber
				\end{eqnarray}
	
				Thus, lemma \ref{5_1_3} tells us that the property is true for those four terms. This step is to be applied repeatedly, to the powers of $g_{\sigma_{n} \sigma_{n-1} ..\sigma_{1} a_{n+1}} $ down to an element of the form $ \tau_{n+1}\big(h (g_{\sigma_{n} \sigma_{n-1} ..\sigma_{1} a_{n+1}})^{1}\big)$, arriving to the sum of:
				
				$$\tau_{n+1}( g_{\sigma_{n} \sigma_{n-1} ..\sigma_{1} a_{n+1}})$$ and $$\tau_{n+1} (h' g_{\sigma_{n-1} ..\sigma_{1} a_{n+1}}) ,$$\\
				
				which is the sum of  values of $ \tau_{n+1}$ on Markov elements, since $h,h'\in F \big(\widehat{TL}_{n} (q)\big)$.
	\end{demo}
			We end this part by the following lemma:			
			
			\begin{lemme} \label{5_1_5}
				Let $ 1\leq k $. Then $ (g_{\sigma_{n} \sigma_{n-1} ..\sigma_{1} a_{n+1}})^{k}$ is a sum of two kinds of elements:\\
				\begin{itemize} 
					\item[$(1)$] $g_{\sigma_{n}} \big(g_{\sigma_{n-1} \sigma_{n-2} .. \sigma_{1}}  F(t_{a_{n}}) \big)^{j} g_{\sigma_{n}} h $, with $j\leq k $.
					\item[$(2)$] $\big(g_{\sigma_{n-1} \sigma_{n-2} .. \sigma_{1}}  F(t_{a_{n}})\big)^{i} g_{\sigma_{n}} f $, with $i < k $,\\
				\end{itemize}	
					
				with $h,f$ in $F \big(\widehat{TL}_{n} (q)\big)$ and $2\leq n$.\\
    
				Moreover, in the first type we have one, and only one element, with $j=k$, in which we have: 
				\vspace{-0.25cm}
				\begin{eqnarray}
					h = \prod^{i=k-1}_{i=0} \phi^{i}\big[F(t^{-1}_{a_{n}})\big]. \nonumber
				\end{eqnarray}
			\end{lemme}
      
			\begin{demo}
				Suppose that $k=1$. Then,
				\begin{eqnarray}
					g_{\sigma_{n} \sigma_{n-1} ..\sigma_{1} a_{n+1}} = g_{\sigma_{n}} \big(g_{\sigma_{n-1} \sigma_{n-2} .. \sigma_{1}} F(t_{a_{n}})\big) g_{\sigma_{n}} F\big(t_{a_{n}}\big)^{-1}, \nonumber
				\end{eqnarray}
				
				the property is true.\\ 

				Suppose the property is true for $k-1$, then, with $ 2\leq k$, we have:
				\begin{eqnarray}
					(g_{\sigma_{n} \sigma_{n-1} ..\sigma_{1} a_{n+1}})^{k} = (g_{\sigma_{n} \sigma_{n-1} ..\sigma_{1} a_{n+1}})^{k-1} g_{\sigma_{n} \sigma_{n-1} ..\sigma_{1} a_{n+1}}. \nonumber
				\end{eqnarray}
				
				We apply the property to $(g_{\sigma_{n} \sigma_{n-1} ..\sigma_{1} a_{n+1}})^{k-1}$, which gives two cases:\\
				
				\begin{itemize} 
					\item[$(1)$] $g_{\sigma_{n}} \big(g_{\sigma_{n-1} \sigma_{n-2} .. \sigma_{1}}  F(t_{a_{n}}) \big)^{j'} g_{\sigma_{n}} h  g_{\sigma_{n} \sigma_{n-1} ..\sigma_{1} a_{n+1}}$, with $ j' \leq k-1 $ which is: \textcolor{white}{...................} $g_{\sigma_{n}} \big(g_{\sigma_{n-1} \sigma_{n-2} .. \sigma_{1}}  F(t_{a_{n}}) \big)^{j'} g_{\sigma_{n}} g_{\sigma_{n} \sigma_{n-1} ..\sigma_{1} a_{n+1}} \psi^{-1}\big[h\big]$, which is equal to:\\
					\begin{eqnarray}
						qg_{\sigma_{n}} &\big(& g_{\sigma_{n-1} \sigma_{n-2} .. \sigma_{1}}  F(t_{a_{n}}) \big)^{j'+1} g_{\sigma_{n}} F\big( (t_{a_{n}})^{-1}\big) \psi^{-1}\big[h\big] \nonumber\\\nonumber\\
						&+& (q-1) g_{\sigma_{n}}  g_{\sigma_{n} \sigma_{n-1} ..\sigma_{1} a_{n+1}}\psi^{-1} \bigg[\big(g_{\sigma_{n-1} \sigma_{n-2} .. \sigma_{1}}  F(t_{a_{n}}) \big)^{j'} \bigg]   \psi^{-1}\left[h\right]. \nonumber
					\end{eqnarray}

					Since, $ j'+1 \leq k $, the first term is clear to be of the first type, while the second term is equal to:
					\begin{eqnarray}
						& & \big(q-1\big)q g_{\sigma_{n-1} ..\sigma_{1}} F\big(t_{a_{n}}\big) g_{\sigma_{n}} F\big( (t_{a_{n}})^{-1}\big)\psi^{-1} \big[\big(g_{\sigma_{n-1} \sigma_{n-2} .. \sigma_{1}}  F(t_{a_{n}}) \big)^{j'} \big]   \psi^{-1}\big[h\big] + \nonumber\\\nonumber\\
						& & \big(q-1\big)^{2} g_{\sigma_{n} \sigma_{n-1} ..\sigma_{1} a_{n+1}}\psi^{-1} \big[\big(g_{\sigma_{n-1} \sigma_{n-2} .. \sigma_{1}}  F(t_{a_{n}}) \big)^{j'} \big] \psi^{-1}\big[h\big] .\nonumber
					\end{eqnarray}

					Here, the first term is of the second type (with $i=1 < k $), and the second term is of the first type (with $j=1$).\\
					
					\item[$(2)$] $\big(g_{\sigma_{n-1} \sigma_{n-2} .. \sigma_{1}}  F(t_{a_{n}})\big)^{i'} g_{\sigma_{n}} f g_{\sigma_{n} \sigma_{n-1} ..\sigma_{1} a_{n+1}}$, with $ i' < k-1 $, which is:
					\begin{eqnarray}
						\big(&g&_{\sigma_{n-1} \sigma_{n-2} .. \sigma_{1}}  F(t_{a_{n}}) \big)^{i'} g_{\sigma_{n}}  g_{\sigma_{n} \sigma_{n-1} ..\sigma_{1} a_{n+1}} \psi^{-1} \big[f\big] = \nonumber\\\nonumber\\
						& & q \big(g_{\sigma_{n-1} \sigma_{n-2} .. \sigma_{1}}  F(t_{a_{n}})\big)^{i'+1} g_{\sigma_{n}} F\big((t_{a_{n}})^{-1}\big) \psi^{-1} \big[f\big] + \nonumber\\\nonumber\\
						& & \big(q-1\big) g_{\sigma_{n}} \big(g_{\sigma_{n-1} ..\sigma_{1}}  F(t_{a_{n}})\big) g_{\sigma_{n}}  F\big((t_{a_{n}})^{-1}\big) \psi^{-1}\big[\big(g_{\sigma_{n-1} \sigma_{n-2} .. \sigma_{1}}  F(t_{a_{n}}) \big)^{i'} \big]  \psi^{-1} \big[f\big].  \nonumber\\\nonumber
					\end{eqnarray}
					
					Since $i'+1 <k $, the first term is of the second type, while the second term is of the first type with $j=1$. The lemma is proven. \\

					(By induction over $k$ again, the last formula is easy). 
				\end{itemize}
				
			\end{demo}
			
			\textbf{Part 2}\\
			
			In this part we will consider a given trace $\tau_{3}$ over $\widehat{TL}_{3} (q)$. The aim is to show that $\tau_{3}$ is uniquely defined by its values on Markov elements. consider
			\begin{eqnarray}
				F_{2}:\widehat{TL}_{2} (q) &\longrightarrow& \widehat{TL}_{3} (q)\nonumber\\
				t_{\sigma_{1}} &\longmapsto& g_{\sigma_{1}}\nonumber\\
				t_{a_{n}} &\longmapsto& g_{\sigma_{2}}g_{a_{3}} g^{-1}_{\sigma_{2}}. \nonumber
			\end{eqnarray}
			
			In this part we will denote   $F_{2}$ by $F$.\\
			
			Lemma \ref{5_1_4} tells that we can uniquely determine $\tau_{3}$ by its values over  $(g_{\sigma_{2}\sigma_{1}a_{3}})^{k} $ for a positive $k$ beside its values on Markov elements. We know as well by lemma \ref{5_1_5} that $(g_{\sigma_{2}\sigma_{1}a_{3}})^{k} $ is a sum of two kinds of elements: \\
			\begin{itemize}[label=$\bullet$, font=\normalsize, font=\color{black}, leftmargin=2cm, parsep=0cm, itemsep=0.25cm, topsep=0cm]
				\item[$(1)$] $g_{\sigma_{2}} \big(g_{\sigma_{1}}  F(t_{a_{2}}) \big)^{j} g_{\sigma_{2}} h $ with $j\leq k $.
				\item[$(2)$] $\big(g_{\sigma_{1}}  F(t_{a_{2}}) \big)^{i} g_{\sigma_{2}} f $ with $i < k $. \\	
			\end{itemize}
			
			Here, $h$ and $f$ are in $F\big(\widehat{TL}_{2}(q)\big)$ .\\
	
			Moreover, in first type, only when $j=k$, we have: 
			\begin{eqnarray}
				h= \prod^{i=k-1}_{i=0} \psi^{i} \bigg[\big(F(t_{a_{2}})^{-1}\big)\bigg]. \nonumber
			\end{eqnarray}
			
			In other terms:
			\begin{eqnarray}
				 g_{\sigma_{2}} \bigg(g_{\sigma_{1}}  F\big(t_{a_{2}}\big) \bigg)^{k} g_{\sigma_{2}}  &=& \bigg(g_{\sigma_{2}\sigma_{1}a_{3}}\bigg)^{k}\prod^{i=k-1}_{i=0} \psi^{i} \bigg[\big(F(t_{a_{2}})\big)\bigg] \nonumber\\\nonumber\\
				 & & - \sum^{r=k-1}_{r=1} \bigg(g_{\sigma_{1}}  F\big(t_{a_{2}}\big) \bigg)^{r} g_{\sigma_{2}} f_{r}\prod^{i=k-1}_{i=0} \psi^{i} \bigg[\big(F(t_{a_{2}})\big)\bigg] \nonumber\\\nonumber\\
				 & & +  \sum^{l=k-1}_{l=1} g_{\sigma_{2}}\bigg(g_{\sigma_{1}}  F\big(t_{a_{2}}\big) \bigg)^{l} g_{\sigma_{2}} f'_{l}\prod^{i=k-1}_{i=0} \psi^{i} \bigg[\big(F(t_{a_{2}})\big)\bigg]. \nonumber\\\nonumber		
			\end{eqnarray}
			
			We repeat the same step on $g_{\sigma_{2}} \big(g_{\sigma_{1}} F(t_{a_{2}})\big)^{l} g_{\sigma_{2}}$ for every $l$. We deduce the following:\\
	
			\begin{corollaire}\label{5_1_6}
			
				For every $h > 0$, we have: $g_{\sigma_{2}} \big(g_{\sigma_{1}}  F(t_{a_{2}})\big)^{h}g_{\sigma_{2}} = \sum\limits^{j=h}_{j=0} c_{j} \big(g_{\sigma_{2}\sigma_{1}a_{3}}\big)^{j} + \sum\limits_{i} M _{i}$. \\ 
				
				Here, $c_{j}$ is in $F\big(\widehat{TL}_{2}(q)\big)$ for every $j$, and $M_{i} $ is a Markov element for every $i$.
			\end{corollaire} 
	
			Our way to prove Theorem \ref{5_1_1} for $n=3$, is to show that $\tau_{3} \big( (g_{\sigma_{2}\sigma_{1}a_{3}})^{k}\big)$ is a linear combination of values of $ \tau_{3} $ on Markov elements and values on elements of the form $c (g_{\sigma_{2}\sigma_{1}a_{3}})^{h}$, where $ h < k $ and $c$ in $F\big(\widehat{TL}_{2}(q)\big)$. Then, using the induction in the proof of Lemma \ref{5_1_4}, beside the fact that $\tau_{3} (g_{\sigma_{2}\sigma_{1}a_{3}}) $ is a linear combination of some values of $\tau_{3}$ on Markov elements, we see that the work is done.\\
		
			\begin{lemme} \label{5_1_7}
				Suppose that $r$ and $s$ are positives ,such that $r \leq s $. Then:
				\begin{eqnarray}
					\tau_{3} \bigg(\big(g_{\sigma_{1}} F(t_{a_{2}})\big)^{r}g_{\sigma_{1}}g_{\sigma_{2}} \big(g_{\sigma_{1}}  F(t_{a_{2}})\big)^{s}g_{\sigma_{1}}g_{\sigma_{2}}\bigg) = \sum^{j=h}_{j=0} c_{j} \bigg(g_{\sigma_{2}\sigma_{1}a_{3}}\bigg)^{j} + \sum_{i} M _{i} ,\nonumber
				\end{eqnarray}
				
				where $h \leq s$, $c_{j}$ is in $F\big(\widehat{TL}_{2}(q)\big)$ for every $j$ and $M_{i} $ is a Markov element for every $i$.
			\end{lemme}
  
			\begin{demo}	
				\begin{eqnarray}
					\tau_{3} \bigg(\big(g_{\sigma_{1}} &F& (t_{a_{2}})\big)^{r}g_{\sigma_{1}}g_{\sigma_{2}} \big(g_{\sigma_{1}}  F(t_{a_{2}})\big)^{s}g_{\sigma_{1}}g_{\sigma_{2}}\bigg) = \nonumber\\\nonumber\\
					&=& \tau_{3} \bigg(\big(g_{\sigma_{1}}  F(t_{a_{2}})\big)^{r}g_{\sigma_{1}}g_{\sigma_{2}} \big(g_{\sigma_{1}} F(t_{a_{2}})\big)^{s}g_{\sigma_{1}} \underbrace{ g_{\sigma_{2} \sigma_{1}a_{3}} \big(g_{\sigma_{2}\sigma_{1}a_{3}}\big)^{-1}} g_{\sigma_{2}}\bigg) \nonumber\\\nonumber\\
					&=& \tau_{3} \bigg( \big(g_{\sigma_{1}}  F(t_{a_{2}})\big)^{r}g_{\sigma_{1}}g_{\sigma_{2}} g_{\sigma_{2}\sigma_{1}a_{3}} \psi \big[ \big(g_{\sigma_{1}}  F(t_{a_{2}})\big)^{s}g_{\sigma_{1}} \big]\big(g_{\sigma_{2}\sigma_{1}a_{3}}\big)^{-1} g_{\sigma_{2}}\bigg) \nonumber\\\nonumber\\
					&=& \tau_{3} \bigg( \big(g_{\sigma_{1}}  F(t_{a_{2}})\big)^{r}g_{\sigma_{1}}g_{\sigma_{2}} g_{\sigma_{2}\sigma_{1}a_{3}} \psi \big[\big(g_{\sigma_{1}}  F(t_{a_{2}})\big)^{s}g_{\sigma_{1}} \big]g^{-1}_{a_{3}}g^{-1}_{\sigma_{1}}\bigg) \nonumber\\\nonumber\\
					&=& \tau_{3} \bigg( \big(g_{\sigma_{1}}  F(t_{a_{2}})\big)^{r-1}g_{\sigma_{1}}g^{2}_{\sigma_{2}}g_{\sigma_{1}}g_{a_{3}} \big( F(t_{a_{2}}) g_{\sigma_{1}}\big)^{s}F\big(t_{a_{2}}\big) g^{-1}_{a_{3}} F\big(t_{a_{2}}\big)\bigg) \nonumber\\\nonumber\\
					&=& \frac{1-q}{q} \tau_{3} \bigg( \big(g_{\sigma_{1}}  F(t_{a_{2}})\big)^{r-1}g_{\sigma_{1}}g^{2}_{\sigma_{2}}g_{\sigma_{1}}g_{a_{3}} \big( F(t_{a_{2}}) g_{\sigma_{1}}\big)^{s}F\big(t_{a_{2}}\big) F\big(t_{a_{2}}\big)\bigg) \nonumber\\\nonumber\\
					&+& \frac{1}{q} \tau_{3} \bigg( \big(g_{\sigma_{1}}  F(t_{a_{2}})\big)^{r-1}g_{\sigma_{1}}g^{2}_{\sigma_{2}}g_{\sigma_{1}}g_{a_{3}} \big( F(t_{a_{2}}) g_{\sigma_{1}}\big)^{s}F\big(t_{a_{2}}\big) g_{a_{3}} F\big(t_{a_{2}}\big)\bigg) . \nonumber			
				\end{eqnarray}
	  
				Now, the term corresponding to $\frac{1-q}{q}$ is $\tau_{3}$  evaluated on the sum of Markov element and an element of style $ c_{j} (g_{\sigma_{2}\sigma_{1}a_{3}})^{1} $. So, We are reduced to the second term, thus, reduced to:
				\begin{eqnarray}
					\tau_{3} \bigg( \big(g_{\sigma_{1}}  F(t_{a_{2}}) \big)^{r-1}g_{\sigma_{1}}g^{2}_{\sigma_{2}}g_{\sigma_{1}}g_{a_{3}} \big( F(t_{a_{2}}) g_{\sigma_{1}}\big)^{s}F \big(t_{a_{2}}\big) g_{a_{3}} F\big(t_{a_{2}}\big)\bigg). \nonumber
				\end{eqnarray}
			
				Obviously, we are in the case:
				\begin{eqnarray}
					q \tau_{3} \bigg( \big( g_{\sigma_{1}}  F(t_{a_{2}})\big)^{r-1} g_{\sigma_{1}}g_{\sigma_{1}}&g&_{a_{3}} \big(F(t_{a_{2}}) g_{\sigma_{1}}\big)^{s}F\big(t_{a_{2}}\big) g_{a_{3}} F\big(t_{a_{2}}\big)\bigg) = \nonumber\\\nonumber\\
					& & q \tau_{3} \bigg( \big(g_{\sigma_{1}}  F(t_{a_{2}})\big)^{r-1}g^{2}_{\sigma_{1}}g_{a_{3}} \big( F(t_{a_{2}}) g_{\sigma_{1}}\big)^{s} F\big(t^{2}_{a_{2}}\big)g_{\sigma_{2}} \bigg),\nonumber
				\end{eqnarray}
				
				since $g_{a_{3}} F\big(t_{a_{2}}\big)= F\big(t_{a_{2}}\big)g_{\sigma_{2}}$.\\
				 
				Now,
				\begin{eqnarray}
					 \tau_{3} \bigg( \big( g_{\sigma_{1}}  F(t_{a_{2}})&\big)&^{r-1}g^{2}_{\sigma_{1}}g_{a_{3}} \big( F(t_{a_{2}}) g_{\sigma_{1}}\big)^{s} F\big(t^{2}_{a_{2}}\big)g_{\sigma_{2}} \bigg)= \nonumber\\\nonumber\\
					& & \big(q-1\big)\tau_{3} \bigg( \big(g_{\sigma_{1}}  F(t_{a_{2}})\big)^{r-1}g^{2}_{\sigma_{1}}g_{a_{3}} F\big(t_{a_{2}}\big) \big(g_{\sigma_{1}} F(t_{a_{2}})\big)^{s} g_{\sigma_{2}} \bigg)\nonumber\\\nonumber\\
					+& & q\tau_{3} \bigg( \big(g_{\sigma_{1}}  F(t_{a_{2}})\big)^{r-1}g^{2}_{\sigma_{1}}g_{a_{3}} F(t_{a_{2}}) g_{\sigma_{1}} \big( F(t_{a_{2}}) g_{\sigma_{1}}\big)^{s-1} g_{\sigma_{2}} \bigg).\nonumber\\\nonumber	
				\end{eqnarray}
	            \vspace{-0.3 cm}
				Which is equal to the sum:
				
				\begin{eqnarray}
					\big(q-1\big)\tau_{3} &\bigg(& \big(g_{\sigma_{1}} F(t_{a_{2}})\big)^{r-1}g^{2}_{\sigma_{1}} F\big(t_{a_{2}}\big) \underbrace{g_{\sigma_{2}} \big( g_{\sigma_{1}} F(t_{a_{2}})\big)^{s} g_{\sigma_{2}}}_{} \bigg) \nonumber\\\nonumber\\
					& & + q\tau_{3} \bigg( \big(g_{\sigma_{1}}  F(t_{a_{2}})\big)^{r-1}g^{2}_{\sigma_{1}} F\big(t_{a_{2}}\big) g_{\sigma_{2}} g_{\sigma_{1}} \big( F(t_{a_{2}}) g_{\sigma_{1}}\big)^{s-1} g_{\sigma_{2}} \bigg).\nonumber
				\end{eqnarray}
				
				Now, the first term is covered by corollary \ref{5_1_6}. Thus we are interested with the second term:
				\begin{eqnarray}
					\tau_{3} \bigg( \big(g_{\sigma_{1}} F(t_{a_{2}})\big)^{r-1}g^{2}_{\sigma_{1}} &F& \big(t_{a_{2}}\big) g_{\sigma_{2}} \big(g_{\sigma_{1}} F(t_{a_{2}})\big)^{s-1} g_{\sigma_{1}} g_{\sigma_{2}} \bigg) = \nonumber\\\nonumber\\
					& & q\tau_{3} \bigg(\big(g_{\sigma_{1}} F(t_{a_{2}})\big)^{r} g_{\sigma_{2}} \big(g_{\sigma_{1}} F(t_{a_{2}})\big)^{s-1} g_{\sigma_{1}} g_{\sigma_{2}} \bigg) \nonumber\\\nonumber\\
					&+& \big(q-1\big)\tau_{3} \bigg(\big(g_{\sigma_{1}} F(t_{a_{2}})\big)^{r-2} g_{\sigma_{1}}F(t^{2}_{a_{2}}) g_{\sigma_{2}} \big(g_{\sigma_{1}}  F(t_{a_{2}})\big)^{s-1} g_{\sigma_{1}} g_{\sigma_{2}}\bigg). \nonumber				
				\end{eqnarray}
				
				Which is equal to:
				\begin{eqnarray}
					q\tau_{3} \bigg( \underbrace{g_{\sigma_{2}} \big( g_{\sigma_{1}} F(t_{a_{2}})\big)^{r}  g_{\sigma_{2}}}_{} &\big(& g_{\sigma_{1}} F(t_{a_{2}})\big)^{s-1} g_{\sigma_{1}} \bigg) + \nonumber\\\nonumber\\
					& & \big(q-1\big) \tau_{3} \bigg(\big(g_{\sigma_{1}} F(t_{a_{2}})\big)^{r-2} g_{\sigma_{1}}F\big(t^{2}_{a_{2}}\big) g_{\sigma_{2}} \big(g_{\sigma_{1}} F(t_{a_{2}})\big)^{s-1} g_{\sigma_{1}} g_{\sigma_{2}} \bigg) \nonumber			
				\end{eqnarray}
				
				The first term is covered by corollary \ref{5_1_6}. We are reduced to
				\begin{eqnarray}
					\bigg(\tau_{3} \big(g_{\sigma_{1}}  F(t_{a_{2}})\big)^{r-2} g_{\sigma_{1}} F \big(t^{2}_{a_{2}}\big) g_{\sigma_{2}} \big(g_{\sigma_{1}}  F(t_{a_{2}})\big)^{s-1} g_{\sigma_{1}} g_{\sigma_{2}} \bigg), \nonumber
				\end{eqnarray}
				
				which is equal to:
				\begin{eqnarray}
					\big(q-1\big)\tau_{3} &\bigg(&\underbrace{g_{\sigma_{2}}\big(g_{\sigma_{1}}  F(t_{a_{2}})\big)^{r-1} g_{\sigma_{2}}}_{} \big(g_{\sigma_{1}} F(t_{a_{2}}) \big)^{s-1} g_{\sigma_{1}}\bigg) + \nonumber\\\nonumber\\
					& & q\tau_{3} \bigg(\big(g_{\sigma_{1}} F(t_{a_{2}})\big)^{r-2} g_{\sigma_{1}} g_{\sigma_{2}} \big(g_{\sigma_{1}} F(t_{a_{2}})\big)^{s-1} g_{\sigma_{1}} g_{\sigma_{2}} \bigg). \nonumber
				\end{eqnarray}
						
				The first term is covered by corollary \ref{5_1_6}. Thus, we see that, in general, the value of $\tau_{3}$ over $\big(g_{\sigma_{1}}  F(t_{a_{2}})\big)^{r} g_{\sigma_{1}} g_{\sigma_{2}} \big(g_{\sigma_{1}}  F(t_{a_{2}}) \big)^{s} g_{\sigma_{1}} g_{\sigma_{2}} $ can be shifted to its value over:
				\begin{eqnarray}
					\big(g_{\sigma_{1}}  F(t_{a_{2}})\big)^{r-2} g_{\sigma_{1}} g_{\sigma_{2}} \big(g_{\sigma_{1}}  F(t_{a_{2}}) \big)^{s-1} g_{\sigma_{1}} g_{\sigma_{2}}. \nonumber
				\end{eqnarray}
				
				After a finite number of repetitions of the computation above (with the possibility of exchanging $r$ and $s$), we see that the lemma is proven modulo determining:
				\begin{eqnarray}
					\tau_{3} \bigg( \big(g_{\sigma_{1}}  F(t_{a_{2}})\big)^{m} \underbrace{g_{\sigma_{1}} g_{\sigma_{2}} g_{\sigma_{1}}}_{} F\big(t_{a_{2}}\big) g_{\sigma_{1}} g_{\sigma_{2}}\bigg). \nonumber
				\end{eqnarray}
				
				We see that the terms corresponding to $-g_{\sigma_{1}}$ and -1 correspond to Markov elements. While those who correspond to $-g_{\sigma_{1}}g_{\sigma_{2}}$ and $-g_{\sigma_{2}}$ are covered by corollary \ref{5_1_6} for $h=1$. Finally the term corresponding to $-g_{\sigma_{2}}g_{\sigma_{1}}$ is covered by corollary \ref{5_1_6} for $h=m$ .
		
			\end{demo}
		
			\begin{lemme} \label{5_1_8}
				Suppose that $r$ and $s$ are positive such that $r \leq s $. Then:
				\begin{eqnarray}
					\tau_{3} \bigg( g_{a_{3}}F\big(t_{a_{2}}\big)\big(g_{\sigma_{1}}  F(t_{a_{2}})\big)^{s}g_{a_{3}} \big(g_{\sigma_{1}}  F(t_{a_{2}})\big)^{r}\bigg) = \sum^{j=h}_{j=0} c_{j} \big(g_{\sigma_{2}\sigma_{1}a_{3}}\big)^{j} + \sum_{i} M _{i}. \nonumber
				\end{eqnarray}
				
				Where $h \leq s$ , $c_{j}$ is in $F(\widehat{TL}_{2}(q))$ for every $j$ and $M_{i} $ is a Markov element for every $i$.
			\end{lemme}
			
			\begin{demo}
				\begin{eqnarray}
					\tau_{3} \bigg( g_{a_{3}} F \big(t_{a_{2}} &\big)& \big(g_{\sigma_{1}}  F(t_{a_{2}})\big)^{s}g_{a_{3}} \big(g_{\sigma_{1}}  F(t_{a_{2}})\big)^{r}\bigg) = \nonumber\\\nonumber\\
					& & \tau_{3} \bigg( g_{a_{3}} \big(g_{\sigma_{2}\sigma_{1}a_{3}}\big)^{-1} g_{\sigma_{2}\sigma_{1}a_{3}}F\big(t_{a_{2}}\big)\big(g_{\sigma_{1}}  F(t_{a_{2}})\big)^{s} g_{a_{3}} \big(g_{\sigma_{1}}  F(t_{a_{2}})\big)^{r}\bigg) = \nonumber\\\nonumber\\
					& & \tau_{3} \bigg( g_{a_{3}}\big(g_{\sigma_{2}\sigma_{1}a_{3}}\big)^{-1} \psi \big[F\big(t_{a_{2}}\big)\big(g_{\sigma_{1}}  F(t_{a_{2}})\big)^{s}\big] g_{\sigma_{2}\sigma_{1}a_{3}}g_{a_{3}} \big(g_{\sigma_{1}}  F(t_{a_{2}})\big)^{r} \bigg) = \nonumber\\\nonumber\\
					& & \tau_{3} \bigg( g^{-1}_{\sigma_{1}} g^{-1}_{\sigma_{2}} g_{\sigma_{1}} \big( F(t_{a_{2}}) g_{\sigma_{1}} \big)^{s} g_{\sigma_{2}}g_{\sigma_{1}}g^{2}_{a_{3}} \big(g_{\sigma_{1}}  F(t_{a_{2}})\big)^{r} \bigg). \nonumber
				\end{eqnarray}

				Here, we see that this term is a sum of two terms coming from  $g^{2}_{a_{3}} = (q-1) g_{a_{3}} + q $. The term corresponding to $ (q-1) g_{a_{3}} $ is covered the same way as in the last lemma (with $a_{3}$ instead of $\sigma_{2}$ above. Hence we treat the term corresponding to $q$, that is: 
				\begin{eqnarray}
					\tau_{3} \bigg( g^{-1}_{\sigma_{1}} \underbrace{g^{-1}_{\sigma_{2}}}_{} \big(g_{\sigma_{1}}  F(t_{a_{2}}) \big)^{s} \underbrace{g_{\sigma_{1}} g_{\sigma_{2}}g_{\sigma_{1}}} \big(g_{\sigma_{1}} F(t_{a_{2}})\big)^{r}\bigg). \nonumber  
				\end{eqnarray}
				
				Before applying TL relations, we see in the same way as above, that we are reduced to the next value (otherwise, it is $\tau_{3}$ evaluated on a Markov element):
				\begin{eqnarray}
					\tau_{3} \bigg( g^{-1}_{\sigma_{1}} \underbrace{g_{\sigma_{2}}}_{} \big(g_{\sigma_{1}}  F(t_{a_{2}})\big)^{s}  \underbrace{g_{\sigma_{1}} g_{\sigma_{2}}g_{\sigma_{1}}} \big(g_{\sigma_{1}}  F(t_{a_{2}})\big)^{r}\bigg). \nonumber
				\end{eqnarray}

				We see that the terms corresponding to $-g_{\sigma_{1}}$ and -1 correspond to Markov elements. And those who correspond to $-g_{\sigma_{2}g_{\sigma_{1}}}$ and $-g_{\sigma_{2}}$, are covered by corollary \ref{5_1_6} for $h=s$. \\
	
				The term corresponding to $-g_{\sigma_{1}}g_{\sigma_{2}}$ is:
				\begin{eqnarray}
					 \tau_{3} \bigg( g^{-1}_{\sigma_{1}} g_{\sigma_{2}} \big(g_{\sigma_{1}}  F(t_{a_{2}})\big)^{s} g_{\sigma_{1}} g_{\sigma_{2}} \big(g_{\sigma_{1}} F(t_{a_{2}})\big)^{r}\bigg), \nonumber
				\end{eqnarray}
				
				which is:
				\begin{eqnarray}
					\frac{1-q}{q}\tau_{3} \bigg( \big(g_{\sigma_{1}} F(t_{a_{2}})\big)^{s} g_{\sigma_{1}} \underbrace{g_{\sigma_{2}} \big(g_{\sigma_{1}}  F(t_{a_{2}})\big)^{r}g_{\sigma_{2}}}_{}\bigg) + \frac{1}{q}\tau_{3} \bigg( g_{\sigma_{1}} g_{\sigma_{2}} \big(g_{\sigma_{1}}  F(t_{a_{2}}) \big)^{s} g_{\sigma_{1}} g_{\sigma_{2}} \big(g_{\sigma_{1}}  F(t_{a_{2}})\big)^{r}\bigg) \nonumber
				\end{eqnarray}
				
				The first term is covered by corollary \ref{5_1_6} for $h=r$. The second follows by lemma \ref{5_1_7}.\\
				
			\end{demo}
	
			Let us go back to $\tau_{3}( g_{\sigma_{2} \sigma_{1} a_{n+1}})^{k}$. The aim is to show that:
			\begin{eqnarray}
				\tau_{3} \big( g_{\sigma_{2} \sigma_{1} a_{n+1}}\big)^{k} =  \tau_{3}\bigg(\sum^{j=h}_{j=0} c_{j} (g_{\sigma_{2}\sigma_{1}a_{3}})^{j} + \Sigma_{i} M _{i} \bigg) ,\nonumber
			\end{eqnarray}
						
			where $h < k $. By lemma \ref{5_1_5}, it is sufficient to deal with:
			\begin{eqnarray}
				g_{\sigma_{2}} \bigg( g_{\sigma_{1}}  F(t_{a_{2}}) \bigg)^{k}g_{\sigma_{2}}\prod^{i=k-1}_{i=0} \psi^{i} \bigg[\big(F(t_{a_{2}})^{-1}\big)\bigg]. \nonumber
			\end{eqnarray}
		
			It is clear that this element is written as a linear combination of four kind of elements:\\
			
			\begin{itemize}[label=$\bullet$, font=\normalsize, font=\color{black}, leftmargin=2cm, parsep=0cm, itemsep=0.25cm, topsep=0cm]
				\item[$1)$] $g_{\sigma_{2}} \big(g_{\sigma_{1}} F(t_{a_{2}}) \big)^{k} g_{\sigma_{2}} \big(g_{\sigma_{1}} F(t_{a_{2}}) \big)^{h}$.
				\item[$2)$] $g_{\sigma_{2}} \big(g_{\sigma_{1}} F(t_{a_{2}}) \big)^{k} g_{\sigma_{2}} \big(g_{\sigma_{1}} F(t_{a_{2}}) \big)^{h}g_{\sigma_{1}}$.
				\item[$3)$] $g_{\sigma_{2}} \big(g_{\sigma_{1}} F(t_{a_{2}}) \big)^{k} g_{\sigma_{2}} \big(  F(t_{a_{2}})g_{\sigma_{1}} \big)^{h} F \big(t_{a_{2}}\big)$.
				\item[$4)$] $g_{\sigma_{2}} \big(g_{\sigma_{1}} F(t_{a_{2}}) \big)^{k} g_{\sigma_{2}}\big(F(t_{a_{2}})g_{\sigma_{1}} \big)^{h} $,  \\			
			\end{itemize}
			
			where $h \leq [\frac{k}{2}] < k $, since $ 1 < k $.\\ 
						
			\begin{itemize}[label=$\bullet$, font=\normalsize, font=\color{black}, leftmargin=2cm, parsep=0cm, itemsep=0.25cm, topsep=0cm]
				\item[$1)$] We start by $\tau_{3} \bigg( g_{\sigma_{2}} \big(g_{\sigma_{1}}  F(t_{a_{2}}) \big)^{k} g_{\sigma_{2}}\big(g_{\sigma_{1}} F(t_{a_{2}}) \big)^{h} \bigg)$. Which is equal to:
					
					\begin{eqnarray}
						\tau_{3} \bigg( \big(g_{\sigma_{1}} F(t_{a_{2}}) \big)^{k} \underbrace{g_{\sigma_{2}} \big(g_{\sigma_{1}} F(t_{a_{2}}) \big)^{h}g_{\sigma_{2}}}_{} \bigg), \nonumber
					\end{eqnarray}
					
					follows directly, regarding corollary \ref{5_1_6} .\\ 
	
				\item[$2)$] Now we consider 				
				\begin{equation}
					\tau_{3} \bigg(g_{\sigma_{2}} \big(g_{\sigma_{1}}  F(t_{a_{2}}) \big)^{k} g_{\sigma_{2}} \big(g_{\sigma_{1}} F(t_{a_{2}}) \big)^{h}g_{\sigma_{1}} \bigg),\label{eq_1} 
				\end{equation}
				
				which is equal to:
					\begin{eqnarray}
						& & \tau_{3} \bigg(g_{\sigma_{2}} \big(g_{\sigma_{1}} F(t_{a_{2}}) \big)^{k} g_{\sigma_{2}} \big(g_{\sigma_{1}}  F(t_{a_{2}}) \big)^{h}g_{\sigma_{1}} g_{\sigma_{2} \sigma_{1} a_{n+1}} \big(g_{\sigma_{2} \sigma_{1} a_{3}} \big)^{-1} \bigg) = \nonumber\\\nonumber\\ 
						& & \tau_{3} \bigg(g_{\sigma_{2}} \big(g_{\sigma_{1}}  F(t_{a_{2}}) \big)^{k} g^{2}_{\sigma_{2}} g_{\sigma_{1}} g_{a_{3}} \psi \big[\big(g_{\sigma_{1}}  F(t_{a_{2}}) \big)^{h}g_{\sigma_{1}}\big]  \big(g_{\sigma_{2} \sigma_{1} a_{3}}\big)^{-1} \bigg)= \nonumber\\\nonumber\\
						& & \tau_{3} \bigg( \big(g_{\sigma_{1}}  F(t_{a_{2}}) \big)^{k-1} g^{2}_{\sigma_{2}} g_{\sigma_{1}} g_{a_{3}} \big(F(t_{a_{2}})g_{\sigma_{1}} \big)^{h}F(t_{a_{2}}) g^{-1}_{a_{3}} F\big(t_{a_{2}}\big)\bigg), \nonumber\\\nonumber				
					\end{eqnarray}
					
					with the very same steps as used above, we see that we are reduced to :	
					\begin{eqnarray}
						& & \tau_{3} \bigg( \big(g_{\sigma_{1}}  F(t_{a_{2}}) )^{k-1} g_{\sigma_{1}} g_{a_{3}} \big( F(t_{a_{2}})g_{\sigma_{1}} \big)^{h}F\big(t_{a_{2}}\big)  g_{a_{3}} F\big(t_{a_{2}}\big) \bigg)\text{, which is: } \nonumber\\\nonumber\\
						& & \tau_{3} \bigg( \big(g_{\sigma_{1}}  F(t_{a_{2}}) \big)^{k-1} g_{\sigma_{1}} \underbrace{g_{a_{3}} F\big(t_{a_{2}}\big)}_{}\big( g_{\sigma_{1}} F(t_{a_{2}}) \big)^{h} \underbrace{g_{a_{3}} F\big(t_{a_{2}}\big)}_{}\bigg)= \nonumber				
					\end{eqnarray}
					
					\begin{eqnarray}
						\tau_{3} \bigg( \big(g_{\sigma_{1}} F(t_{a_{2}}) \big)^{k-1} g_{\sigma_{1}} \underbrace{F\big(t_{a_{2}}\big) g_{\sigma_{2}}}_{} \big( g_{\sigma_{1}} F(t_{a_{2}}) \big)^{h} \underbrace{ F\big(t_{a_{2}}\big) g_{\sigma_{2}}}_{}\bigg), \text{ which is equal to:} \nonumber
					\end{eqnarray}
					\vspace{-0.5cm}
					\begin{footnotesize}
						\begin{eqnarray}
							(q-1 &)&\tau_{3} \bigg( \big( g_{\sigma_{1}}  F(t_{a_{2}}) \big)^{k} g_{\sigma_{2}} \big( g_{\sigma_{1}} F(t_{a_{2}}) \big)^{h} g_{\sigma_{2}}\bigg) + q\tau_{3} \bigg( \big(g_{\sigma_{1}} F(t_{a_{2}}) \big)^{k} g_{\sigma_{2}} \big( g_{\sigma_{1}} F(t_{a_{2}}) \big)^{h-1}  g_{\sigma_{1}} g_{\sigma_{2}}\bigg),\nonumber
						\end{eqnarray}
					\end{footnotesize}
		
					we see that corollary \ref{5_1_6} covers the first term. Thus we see that:
					\begin{eqnarray}
						\tau_{3}\bigg(g_{\sigma_{2}} \big(g_{\sigma_{1}}  F(t_{a_{2}}) \big)^{k} g_{\sigma_{2}} \big(g_{\sigma_{1}} F(t_{a_{2}}) \big)^{h}g_{\sigma_{1}}\bigg)\text{ in (eq. \ref{eq_1}), is shifted to: }\nonumber
					\end{eqnarray}
					\vspace{-0.3 cm}		
					\begin{eqnarray}
						\tau_{3} \bigg(g_{\sigma_{2}} \big(g_{\sigma_{1}}  F(t_{a_{2}}) \big)^{k} g_{\sigma_{2}} \big(g_{\sigma_{1}} F(t_{a_{2}}) \big)^{h-1}g_{\sigma_{1}} \bigg), \nonumber
					\end{eqnarray}
					
					going on in this manner, we arrive to: 
					\begin{eqnarray}
						\tau_{3} \bigg(g_{\sigma_{2}} \big(g_{\sigma_{1}}  F(t_{a_{2}}) \big)^{k} g_{\sigma_{2}}g_{\sigma_{1}} F\big(t_{a_{2}}\big) g_{\sigma_{1}} \bigg), \nonumber
					\end{eqnarray}
					
					with the same steps above, we see that we are reduced to
					\begin{eqnarray}
						\tau_{3} \bigg( \big(g_{\sigma_{1}} F\big(t_{a_{2}}\big) \big)^{k} g_{\sigma_{2}}g_{\sigma_{1}}  F(t^{2}_{a_{2}})g_{\sigma_{2}}\bigg), \text{ which is equal to } \nonumber\\\nonumber\\
						\big(q-1\big)\tau_{3} \bigg( \big(g_{\sigma_{1}}  F(t_{a_{2}}) \big)^{k} g_{\sigma_{2}}g_{\sigma_{1}} F\big(t_{a_{2}}\big)g_{\sigma_{2}}\bigg) + q\tau_{3} \bigg( \big(g_{\sigma_{1}}  F(t_{a_{2}}) \big)^{k} g_{\sigma_{2}}g_{\sigma_{1}} g_{\sigma_{2}} \bigg), \nonumber
					\end{eqnarray} 
					
					corollary \ref{5_1_6} and TL relations end the job.\\
				
				\item[$3)$] Here we deal with $\tau_{3} \bigg(g_{\sigma_{2}} \big(g_{\sigma_{1}}  F(t_{a_{2}}) \big)^{k} g_{\sigma_{2}}  \big( F(t_{a_{2}})g_{\sigma_{1}} \big)^{h}F\big(t_{a_{2}}\big) \bigg)$, which is:
					
					\vspace{-0.5cm}
					
					\begin{eqnarray}	
						\tau_{3} \bigg( F\big(t_{a_{2}} &\big)& g_{\sigma_{2}} \big(g_{\sigma_{1}}  F(t_{a_{2}}) \big)^{k} g_{\sigma_{2}} \big( F(t_{a_{2}})g_{\sigma_{1}} \big)^{h} \bigg) \nonumber\\\nonumber\\
						&=& \tau_{3} \bigg( g_{a_{3}} F\big(t_{a_{2}}\big) \big(g_{\sigma_{1}}  F(t_{a_{2}}) \big)^{k} g_{\sigma_{2}} \big( F(t_{a_{2}})g_{\sigma_{1}} \big)^{h} \bigg)\nonumber\\\nonumber\\
						&=& \tau_{3} \bigg( g_{a_{3}} \big(F(t_{a_{2}})g_{\sigma_{1}} \big)^{k} F\big(t_{a_{2}}\big) g_{\sigma_{2}} \big(F(t_{a_{2}})g_{\sigma_{1}} \big)^{h} \bigg), \nonumber
					\end{eqnarray}
	
	                but, $g_{\sigma_{2}} = F\big(t^{-1}_{a_{2}}\big) g_{a_{3}} F\big(t_{a_{2}}\big)$, thus:
					
					\vspace{0.6 cm}
					
					\begin{eqnarray}
						& & \tau_{3} \bigg( g_{a_{3}} \big(F(t_{a_{2}})g_{\sigma_{1}} \big)^{k} g_{a_{3}}  F\big(t_{a_{2}}\big) \big(F(t_{a_{2}})g_{\sigma_{1}} \big)^{h} \bigg)  \nonumber\\\nonumber\\
						=& & \tau_{3} \bigg( g_{a_{3}}g^{-1}_{\sigma_{2} \sigma_{1} a_{3}}g_{\sigma_{2} \sigma_{1} a_{3}} \big(F(t_{a_{2}})g_{\sigma_{1}} \big)^{k} g_{a_{3}} F\big(t_{a_{2}}\big) \big(F(t_{a_{2}})g_{\sigma_{1}} \big)^{h} \bigg)  \nonumber\\\nonumber\\
						=& & \tau_{3} \bigg( g^{-1}_{\sigma_{2}} \psi \big[ \big(F(t_{a_{2}})g_{\sigma_{1}}   \big)^{k}\big] g_{\sigma_{2} \sigma_{1} a_{3}} g_{a_{3}} F\big(t_{a_{2}}\big) \big( F(t_{a_{2}})g_{\sigma_{1}} \big)^{h} g^{-1}_{\sigma_{1}} \bigg), \nonumber
					\end{eqnarray}
					
					as we have done above, using the quadratic relations, we see that we are reduced to:
					
					\vspace{0.2cm}
					\begin{eqnarray}
						& & \tau_{3}\bigg(g_{\sigma_{2}} \psi \big[ \big(F(t_{a_{2}})g_{\sigma_{1}} \big)^{k} \big] g_{\sigma_{2}} g_{\sigma_{1}} F\big(t_{a_{2}}\big) \big( F(t_{a_{2}})g_{\sigma_{1}} \big)^{h}  g^{-1}_{\sigma_{1}}\bigg)  \nonumber\\\nonumber\\
						=& & \tau_{3} \bigg(g_{\sigma_{2}} \psi \big[ \big(F(t_{a_{2}})g_{\sigma_{1}} \big)^{k} \big] g_{\sigma_{2}} g_{\sigma_{1}} F\big(t^{2}_{a_{2}}\big) \big( g_{\sigma_{1}} F(t_{a_{2}})\big)^{h-1} \bigg)  \nonumber\\\nonumber\\
						=& & \big(q-1\big)\tau_{3} \bigg(g_{\sigma_{2}} \psi \big[ \big(F(t_{a_{2}})g_{\sigma_{1}} \big)^{k} \big] g_{\sigma_{2}} \big( g_{\sigma_{1}}  F(t_{a_{2}}) \big)^{h} \bigg) + \nonumber\\\nonumber\\
						& & q\tau_{3} \bigg( g_{\sigma_{2}} \psi \big[ \big(F(t_{a_{2}})g_{\sigma_{1}} \big)^{k} \big] g_{\sigma_{2}} g_{\sigma_{1}} \big( g_{\sigma_{1}} F(t_{a_{2}}) \big)^{h-1} \bigg), \nonumber
					\end{eqnarray}
					\vspace{0.3cm}
					
					the first term is covered by corollary \ref{5_1_6}. For the second we see that it is equal to:		
					\begin{eqnarray}
						q\tau_{3} \bigg(g_{\sigma_{2}} \big(g_{\sigma_{1}} &F& (t_{a_{2}})\big)^{k} g_{\sigma_{2}} g_{\sigma_{1}} \big( g_{\sigma_{1}}  F(t_{a_{2}}) \big)^{h-1} \bigg), \text{ which is equal to} \nonumber\\\nonumber\\
						& & q\big(q-1\big)\tau_{3}\bigg(g_{\sigma_{2}} \big(g_{\sigma_{1}} F(t_{a_{2}}) \big)^{k} g_{\sigma_{2}} \big( g_{\sigma_{1}}  F(t_{a_{2}})\big)^{h-1} \bigg) + \nonumber\\\nonumber\\
						& & q^{2}\tau_{3} \bigg(g_{\sigma_{2}} \big(g_{\sigma_{1}} F(t_{a_{2}}) \big)^{k} g_{\sigma_{2}} F\big(t_{a_{2}}\big) \big(g_{\sigma_{1}}  F(t_{a_{2}})\big)^{h-2} \bigg),\nonumber\\\nonumber			
					\end{eqnarray}
					
					the first term is obviously, covered by corollary \ref{5_1_6}, for the second one we see that it is case 3 itself, but with $ h-2 $ instead of $h$. Thus, we get two elements for $\tau_{3}$ to be evaluated on:\\ 
					
					\begin{itemize}[label=$\bullet$, font=\normalsize, font=\color{black}, leftmargin=2cm, parsep=0cm, itemsep=0.25cm, topsep=0cm]
						\item[[$~a$]] $g_{\sigma_{2}} \big( g_{\sigma_{1}} F(t_{a_{2}}) \big)^{k} g_{\sigma_{2}}F\big(t_{a_{2}}\big) g_{\sigma_{1}} F\big(t_{a_{2}}\big)$, 
						\item[[$~b$]] $g_{\sigma_{2}} \big(g_{\sigma_{1}} F(t_{a_{2}}) \big)^{k} g_{\sigma_{2}} F\big(t_{a_{2}}\big) \big( g_{\sigma_{1}} F(t_{a_{2}})\big)^{2}$.\\
					\end{itemize}
					
					For [$~b~$] we can repeat what we have done until arriving to:
					\begin{eqnarray}
						\tau_{3} \bigg(g_{\sigma_{2}} \big( g_{\sigma_{1}} F(t_{a_{2}}) \big)^{k}g_{\sigma_{2}} g_{\sigma_{1}} g_{\sigma_{1}} F\big(t_{a_{2}}\big) \bigg), \text{ which is the following sum:} \nonumber
					\end{eqnarray}
					\vspace{-0.5cm}
					\begin{eqnarray}	
						\big(q-1\big)\tau_{3} \bigg( \big(g_{\sigma_{1}} F(t_{a_{2}}) \big)^{k} g_{\sigma_{2}} g_{\sigma_{1}} F\big(t_{a_{2}}\big) g_{\sigma_{2}} \bigg) + q\tau_{3} \bigg( \big( g_{\sigma_{1}} F(t_{a_{2}}) \big)^{k}\underbrace{g_{\sigma_{2}} F\big(t_{a_{2}} \big) g_{\sigma_{2}}}_{F(t_{a_{2}}) g_{\sigma_{2}}F(t_{a_{2}})} \bigg), \nonumber				
					\end{eqnarray}
					
					obviously, the first term is covered by corollary \ref{5_1_6}, the second term is a Markov element.\\
	
					For [$~a~$] we see that: 
					\begin{eqnarray}
						\tau_{3} \bigg( g_{\sigma_{2}} &\big(& g_{\sigma_{1}} F(t_{a_{2}}) \big)^{k} g_{\sigma_{2}} F\big(t_{a_{2}}\big) g_{\sigma_{1}} F\big(t_{a_{2}}\big) \bigg)  \nonumber\\\nonumber\\
						=\tau_{3} \bigg( g_{\sigma_{2}} &\big(& g_{\sigma_{1}} F(t_{a_{2}})\big)^{k-1}  g_{\sigma_{1}} \underbrace{ F\big(t_{a_{2}}\big) g_{\sigma_{2}} F\big(t_{a_{2}}\big)}_{g^{2}_{\sigma_{2}} g_{a_{3}}} g_{\sigma_{1}} F\big(t_{a_{2}}\big) \bigg)  \nonumber\\\nonumber\\
						=& & \big(q-1\big) \tau_{3} \bigg( g_{\sigma_{2}} \big(g_{\sigma_{1}} F(t_{a_{2}}) \big)^{k-1}  g_{\sigma_{1}} g_{\sigma_{2}} g_{a_{3}} g_{\sigma_{1}} F\big(t_{a_{2}}\big) \bigg) + \nonumber\\\nonumber\\
						& & q\tau_{3} \bigg( g_{\sigma_{2}} \big(g_{\sigma_{1}} F(t_{a_{2}}) \big)^{k-1}  g_{\sigma_{1}}g_{a_{3}} g_{\sigma_{1}} F\big(t_{a_{2}}\big) \bigg),\nonumber		
					\end{eqnarray}
					
					the first term is covered by corollary \ref{5_1_6}, since it is equal to:
					\begin{eqnarray}
						\big(q-1\big) \tau_{3} \bigg( \big( g_{\sigma_{1}} F(t_{a_{2}}) \big)^{k-1} g_{\sigma_{1}} F\big(t_{a_{2}}\big) g_{\sigma_{2}} g_{\sigma_{1}} F\big(t_{a_{2}}\big) g_{\sigma_{2}} \bigg).\nonumber
					\end{eqnarray}

					For the second term, we see that:
					\begin{eqnarray}
						& & q\tau_{3} \bigg( g_{\sigma_{2}} \big(g_{\sigma_{1}} F(t_{a_{2}}) \big)^{k-1}  g_{\sigma_{1}} g_{a_{3}} g_{\sigma_{1}} F\big(t_{a_{2}}\big) \bigg)  \nonumber\\\nonumber\\
						& & =q\tau_{3} \bigg( \big(g_{\sigma_{1}} F(t_{a_{2}}) \big)^{k-1} g_{\sigma_{1}} g_{a_{3}} g_{\sigma_{1}} F\big(t_{a_{2}}\big) g_{\sigma_{2}} \bigg) \nonumber
					\end{eqnarray}
					
					\begin{eqnarray}
						=q\tau_{3} \bigg( \big( g_{\sigma_{1}} F(t_{a_{2}}) \big)^{k-1} g_{\sigma_{1}} g_{a_{3}} g_{\sigma_{1}} g_{a_{3}} F\big(t_{a_{2}}\big) \bigg)  \nonumber\\\nonumber\\
						=q\tau_{3} \bigg( \big( g_{\sigma_{1}} F(t_{a_{2}}) \big)^{k-1} g^{2}_{\sigma_{1}} g_{a_{3}} g_{\sigma_{1}} F\big(t_{a_{2}} \big) \bigg), \nonumber
					\end{eqnarray}
					which is a Markov element, since $ g_{a_{3}} =  F(t_{a_{2}})g_{\sigma_{2}} F(t^{-1}_{a_{2}})$ .\\
					
				\item[$4)$] We deal with  $\tau_{3} \bigg( g_{\sigma_{2}} \big( g_{\sigma_{1}}  F(t_{a_{2}}) \big)^{k} g_{\sigma_{2}} \big( F(t_{a_{2}}) g_{\sigma_{1}} \big)^{h} \bigg) $, using the same techniques:

					\begin{eqnarray}
						& & \tau_{3} \bigg( g_{\sigma_{2}} \big( g_{\sigma_{1}}  F(t_{a_{2}}) \big)^{k} g_{\sigma_{2}} \big( F(t_{a_{2}})g_{\sigma_{1}} \big)^{h} \bigg)  \nonumber\\\nonumber\\
						=& & \tau_{3} \bigg( g_{\sigma_{2}} \big(g_{\sigma_{1}}  F(t_{a_{2}}) \big)^{k} g_{\sigma_{2}} \big( F(t_{a_{2}})g_{\sigma_{1}} \big)^{h} g_{\sigma_{2} \sigma_{1} a_{3}}g^{-1}_{\sigma_{2} \sigma_{1} a_{3}} \bigg)  \nonumber\\\nonumber\\
						=& & \tau_{3} \bigg( g_{\sigma_{2}} \big( g_{\sigma_{1}}  F(t_{a_{2}}) \big)^{k} g^{2}_{\sigma_{2}} g_{\sigma_{1}} g_{a_{3}} \big(g_{\sigma_{1}}F(t_{a_{2}}) \big)^{h} g^{-1}_{\sigma_{2} \sigma_{1} a_{3}} \bigg), \nonumber				
					\end{eqnarray}

					so, we are reduced to:
					\begin{eqnarray}
						& & \tau_{3} \bigg( \big( g_{\sigma_{1}}  F(t_{a_{2}}) \big)^{k-1} g_{\sigma_{1}} g_{a_{3}} \big( g_{\sigma_{1}} F(t_{a_{2}}) \big)^{h} g_{a_{3}} F\big(t_{a_{2}}\big) \bigg). \text{ Which is equal to: } \nonumber\\\nonumber\\
						& & \tau_{3} \bigg( \big( g_{\sigma_{1}}  F(t_{a_{2}}) \big)^{k-1} \underbrace{ g_{\sigma_{1}} g_{a_{3}}g_{\sigma_{1}}}_{V(g_{\sigma_{1}}, g_{a_{3}})} F\big(t_{a_{2}}\big) \big(g_{\sigma_{1}}F(t_{a_{2}}) \big)^{h-1}  F\big(t_{a_{2}}\big) g_{\sigma_{2}} \bigg), \nonumber					
					\end{eqnarray}

					for -1 and $-g_{\sigma_{1}}$ it is a Markov element. For $- g_{a_{3}}g_{\sigma_{1}}$ we see that:

					\begin{eqnarray}
						& & \tau_{3} \bigg( \big( g_{\sigma_{1}} F(t_{a_{2}}) \big)^{k-1} g_{a_{3}} g_{\sigma_{1}} F\big(t_{a_{2}}\big) \big(g_{\sigma_{1}}F(t_{a_{2}}) \big)^{h-1} F\big(t_{a_{2}}\big) g_{\sigma_{2}}\bigg)  \nonumber\\\nonumber\\
						=& & \tau_{3} \bigg( g_{a_{3}} F\big(t_{a_{2}}\big) \big(g_{\sigma_{1}}  F(t_{a_{2}}) \big)^{k-1} g_{a_{3}} \big( g_{\sigma_{1}} F(t_{a_{2}}) \big)^{h} \bigg), \nonumber
					\end{eqnarray}
					
					which is covered by lemma \ref{5_1_8}.\\ 
					
					For $- g_{a_{3}}$, we see that:
					\begin{eqnarray}
						\tau_{3} &\bigg(& \big(g_{\sigma_{1}}  F(t_{a_{2}}) \big)^{k-1} g_{a_{3}} F\big(t_{a_{2}} \big) \big( g_{\sigma_{1}}F(t_{a_{2}}) \big)^{h-1} F\big(t_{a_{2}}\big)g_{\sigma_{2}} \bigg)  \nonumber\\\nonumber\\
						=\tau_{3} &\bigg(& \big( g_{\sigma_{1}} F(t_{a_{2}}) \big)^{k-2} g_{\sigma_{1}}  F\big(t^{2}_{a_{2}}\big) g_{\sigma_{2}} \big(g_{\sigma_{1}}F(t_{a_{2}}) \big)^{h-2} g_{\sigma_{1}} F\big(t^{2}_{a_{2}}\big) g_{\sigma_{2}} \bigg) \nonumber\\\nonumber\\
						 =&\big(&q-1\big) \tau_{3} \bigg( \big(g_{\sigma_{1}}  F(t_{a_{2}}) \big)^{k-1} g_{\sigma_{2}} \big(g_{\sigma_{1}}F(t_{a_{2}}) \big)^{h-2} g_{\sigma_{1}} F\big(t^{2}_{a_{2}} \big) g_{\sigma_{2}} \bigg) + \nonumber\\\nonumber\\
						q\tau_{3} &\bigg(& \big( g_{\sigma_{1}}  F(t_{a_{2}}) \big)^{k-2}g_{\sigma_{1}}   g_{\sigma_{2}} \big(g_{\sigma_{1}}F(t_{a_{2}}) \big)^{h-2} g_{\sigma_{1}} F\big(t^{2}_{a_{2}}\big) g_{\sigma_{2}} \bigg), \nonumber
					\end{eqnarray}
					
					the first term is covered by corollary \ref{5_1_6}. We do the same thing with $F(t^{2}_{a_{2}})$ in the second term, we arrive to:
					\begin{eqnarray}
						q^{2} \tau_{3} \bigg( \big( g_{\sigma_{1}}  F(t_{a_{2}}) \big)^{k-2}g_{\sigma_{1}}   g_{\sigma_{2}} \big( g_{\sigma_{1}}F(t_{a_{2}}) \big)^{h-2} g_{\sigma_{1}} g_{\sigma_{2}} \big), \nonumber
					\end{eqnarray}
					
					which is the case of lemma \ref{5_1_7}.  \\

					For $-g_{\sigma_{1}} g_{a_{3}}$ we see that: 
					\begin{eqnarray}
						\tau_{3} \bigg( \big( g_{\sigma_{1}} F(t_{a_{2}}) &\big)&^{k-1} g_{\sigma_{1}} g_{a_{3}} F\big(t_{a_{2}}\big) \big(g_{\sigma_{1}}F(t_{a_{2}}) \big)^{h-1} F\big(t_{a_{2}}\big) g_{\sigma_{2}} \bigg)  \nonumber\\\nonumber\\ 
						=\tau_{3} \bigg( \big( g_{\sigma_{1}} F(t_{a_{2}}) &\big)&^{k} g_{\sigma_{2}} \big( g_{\sigma_{1}} F(t_{a_{2}}) \big)^{h-1} F\big(t_{a_{2}}\big) g_{\sigma_{2}} \bigg)  \nonumber\\\nonumber\\
						=&\big(& q-1\big) \tau_{3} \bigg( \big( g_{\sigma_{1}} F(t_{a_{2}}) \big)^{k} g_{\sigma_{2}} \big( g_{\sigma_{1}} F(t_{a_{2}}) \big)^{h-1} g_{\sigma_{2}} \bigg) + \nonumber\\\nonumber\\
						q\tau_{3} &\bigg(& \big( g_{\sigma_{1}} F(t_{a_{2}}) \big)^{k} g_{\sigma_{2}} \big(g_{\sigma_{1}}F(t_{a_{2}}) \big)^{h-2} g_{\sigma_{1}} g_{\sigma_{2}} \bigg), \nonumber					
					\end{eqnarray}
					
					corollary \ref{5_1_6} covers the first term, while the second term is covered by (1) from our four cases. \\

			\end{itemize}
			

	\vspace{2cm}

\renewcommand{\refname}{REFERENCES}


\begin{thebibliography}{}

\bibitem{Graham_Lehrer_1998} J. J. Graham and G. I. Lehrer. The representation theory of affine Temperley-Lieb algebras. L'Enseignement Mathematique, 44, 173-218, 1998.\label{Graham_Lehrer_1998}


\bibitem{Graham_Lehrer_2003} J. J. Graham and G. I. Lehrer. Diagram algebras, Hecke algebras and decomposition numbers at roots of unity. Annales Scientifiques de lÉcole Normale Supérieure, 36, Issue 4:479-524, 2003. \label{Graham_Lehrer_2003}

\bibitem{Sadek_2013_1} S. Al Harbat.Markov trace on a tower of affine Temperley-Lieb algebras of type $\tilde{A_{n}}$. 2013. \href{http://arxiv.org/abs/1311.7092}{\emph{ 	arXiv:1311.7092}} \label{Sadek_2013_1}

\bibitem{Sadek_2013_2} S. Al Harbat.A classification of affine fully commutative elements. 2013. \href{http://xxx.tau.ac.il/abs/1311.7089}{\emph{arXiv:1311.7089v1}} \label{Sadek_2013_2}


\bibitem{Sadek_Thesis} S. Al Harbat. On the affine braid group, affine Temperley-Lieb algebra and Markov trace. PH.D
Thesis, 2013. \label{Sadek_Thesis}


\bibitem{Jones_1985} V. F. R. Jones. A polynomial invariant for knots via Von Neumann algebras. Bulletin, American
Mathematical Society, 12, No. 1:103-111, 1985. \label{Jones_1985}


\end{thebibliography}
\end{document}